\documentclass[preprint,11pt]{elsarticle}
\usepackage{epsfig}
\usepackage{subfig}
\pdfoutput=1 
\usepackage{amsmath}
\usepackage{todo}
\usepackage{marginnote}
\usepackage{amsfonts}
\usepackage{amssymb}
\usepackage{graphicx}%
\usepackage[margin=1in]{geometry}
\usepackage{epstopdf}
\providecommand{\U}[1]{\protect\rule{.1in}{.1in}}

\begin{document}
\begin{frontmatter}


\title{On optimal performance of nonlinear energy sinks in multiple-degree-of-freedom systems}


\author{Astitva Tripathi}
\address{Purdue University}
\author {Piyush Grover\corref{PG}}
\address{Mitsubishi Electric Research Labs, Cambridge, MA, USA}
\ead{grover@merl.com}
\author{Tam\'{a}s Kalm\'{a}r-Nagy\corref{TK}} 
\address{Department of Fluid Mechanics, Faculty of
Mechanical Engineering, Budapest University of Technology and Economics,
Hungary}
\cortext[PG]{ Corresponding Author}

\begin{abstract}
We study the problem of optimizing the performance of a nonlinear spring-mass-damper attached to a class of multiple-degree-of-freedom systems. We aim to maximize the rate of one-way energy transfer from primary system to the attachment, and focus on impulsive excitation of a two-degree-of-freedom primary system with an essentially nonlinear attachment. The nonlinear attachment is shown to be able to perform as a `nonlinear energy sink' (NES) by taking away energy from the primary system irreversibly for some types of impulsive excitations. Using perturbation analysis and exploiting separation of time scales, we perform dimensionality reduction of this strongly nonlinear system. Our analysis shows that efficient energy transfer to nonlinear attachment in this system occurs for initial conditions close to homoclinic orbit of the slow time-scale undamped system, a phenomenon that has been previously observed for the case of single-degree-of-freedom primary systems. Analytical formulae for optimal parameters for given impulsive excitation input are derived. Generalization of this framework to systems with arbitrary number of degrees-of-freedom of the primary system is also discussed. The performance of both linear and nonlinear optimally tuned attachments is compared. While NES performance is sensitive to magnitude of the initial impulse, our results show that NES performance is more robust than linear tuned-mass-damper to several parametric perturbations. Hence, our work provides evidence that homoclinic orbits of the underlying Hamiltonian system play a crucial role in efficient nonlinear energy transfers, even in high dimensional systems, and gives new insight into robustness of systems with essential nonlinearity. 
\end{abstract}

\begin{keyword}
nonlinear energy sinks \sep optimization \sep homoclinic orbits \sep multiple time scales \sep resonant capture \sep  averaging theory \sep hamiltonian systems
%
%
\end{keyword}

\end{frontmatter}

\section{Introduction}

The suppression of vibrational energy via transfer from the main structure to
an attachment -both actively or passively- has been a lively research area
since the seminal invention of the tuned mass damper (TMD) \cite{frahm1911device}. With advances in electro-mechanical devices, active
control schemes offer the best performance in terms of vibration absorption.
However, in addition to cost and energy consumption associated with active
control, robustness and stability are crucial factors. Passive vibration
reduction approaches include direct use or variations of linear tuned mass
dampers. However, even if the tuned mass damper is initially tuned to
eliminate resonant responses near the eigenfrequency of the primary system,
the mitigating performance may become less effective due to natural mistuning
of the system parameters (e.g. the varying mass of the secondary due to time
varying load). In H$_{\infty}/$H$_{2}$ optimization the TMD is designed such
that the maximum amplitude magnification factor or the squared area under the
response curve of the primary system is minimized, respectively. Analytical solutions for the
H$_{\infty}/$H$_{2}$ optimization of the TMD have been found \cite{asami2002analytical}, in the form of series solution for
the H$_{\infty}$ optimization and a closed-form algebraic solution for the
H$_{2}$ optimization. In related work \cite{zilletti2012optimisation}, an optimization problem is considered which either minimizes the kinetic energy of the host
structure or maximizes the power dissipation within the absorber. Formulas for the optimal ratio of the absorber natural
frequency to the host natural frequency and optimal damping ratio of the
absorber were also obtained in that work. Ref. \cite{bisegna2012closed} deals with the
analysis and optimization of tuned mass dampers by providing design formulas
for maximizing the exponential time-decay rate of the system transient
response. A detailed analysis is presented for the classical TMD
configuration, involving an auxiliary mass attached to the main structure by
means of a spring and a dashpot. Analytic expressions of the optimal
exponential time-decay rate are obtained for any mass ratio and tuning
condition. Then, a further optimization with respect to the latter is
performed. 

Lyapunov's second method has been used
to minimize an integral square performance measure of damped vibrating
structures subject to initial impulse \cite{wang1984transient}. Using the same approach, the closed-form solutions of optimum
parameters for undamped primary structure utilizing Kronecker product and
matrix column expansion are derived in Ref. \cite{du2008analytical}. Other related works include a
parametric study on a TMD using steady-state harmonic excitation analysis and
time-history analysis (with the El Centro and Mexico earthquake excitation
signals) \cite{rana1998parametric}, H$_{\infty}$
optimal design of a dynamic vibration absorber variant (the damping element is
connected directly to the ground instead of the primary mass) for suppressing
high-amplitude vibrations of damped primary systems \cite{chun2015h}, and 
determination of optimal absorber parameters to maximize of the primary system frequency response \cite{harik2013design}.
A simple method for choosing optimal parameters for a two-degree-of-freedom (translational/rotational) TMD has been reported in \cite{jang2012simple}. This method uses the fixed points of the frequency
response functions to determine the stiffness of the TMD for a given mass. The effectiveness of TMDs in reducing the transient structural response for impulsive loadings has also been investigated in Ref. \cite{salvi2013analysis}.

Nonlinear energy transfer between modes due to resonance has also been studied
extensively \cite{nayfeh2008nonlinear}, focusing on modal interactions and
transfers from high to low frequency modes. The energy transfer phenomenon in
this class of systems is essentially modal, and does not necessarily translate
to one-way transfer between spatially distinct components of the system.
Transition between resonances in Hamiltonian systems has been studied via
geometrical and analytical methods in the past few decades
\cite{koon2000heteroclinic,vainchtein2006passage}. Recently, use of active
control strategies to move the system between these resonances have been
explored \cite{schroer1997targeting, vainchtein2004capture,grover2009designing,grover2012optimized}.

Targeted energy transfers (TETs), i.e. passively controlled transfers of
vibrational energy in coupled oscillators to a targeted component where the
energy eventually localizes, have been a topic of great interest in the past decade
\cite{vakbook}. The basic device is called a nonlinear energy sink (NES),
which generally consists of a light mass, an essentially nonlinear spring and
a viscous damper. Properly designed, the NES is capable of one-way channeling
of unwanted energy from a primary system to NES over broadband frequency
ranges. TET is realized through resonance captures and escapes from
resonances, following (countable infinite) resonance manifolds due to the
essential nonlinearity. While the phenomenon of targeted energy transfer has
been extensively studied in this context
\cite{vakakis2001energy,vakakis2003dynamics,kerschen2005irreversible}, the
parameter selection and optimization problem for multiple-degree-of-freedom
systems is still a challenge. In Ref. \cite{sapsis_01}, energy transfer initiated by an impulsive input in a single-degree-of-freedom system coupled with NES was analyzed and the optimal energy transfer phenomenon was described in terms of existence of a homoclinic orbit in a reduced phase space of the undamped (Hamiltonian) averaged (slow) system.

In this paper, we extend this analysis of Ref. \cite{sapsis_01} to a class of weakly damped multiple-degree-of-freedom systems, especially focusing on a
two-degree-of-freedom system with an attached NES. We obtain the near-optimal parameters for the NES using the complexification-averaging technique and slow
flow analysis. The validity of dimensionality reduction enabled by our analysis is supported by numerical comparisons between original and reduced order systems. Using a combination of perturbation analysis and simulation, we show that under assumptions of weak damping in the linear system, the homoclinic orbit picture persists in higher degree-of-freedom systems. We use Lyapunov analysis to optimize a linear TMD using a similar cost function, and obtain a semi-analytical solution for optimal parameters. Using the semi-analytical formulae, we are able to perform extensive performance comparison studies using these two classes of optimally tuned vibration absorbers.

The structure of the paper is as follows. In Section \ref{two-dof}, we consider a two-degree-of-freedom linear system (called primary) with an attached NES. We perform a numerical study to compute
various branches of periodic solutions, and obtain the frequency-energy plot.
Focusing on the 1:1:1 resonance between the two masses of the primary system,
and the nonlinear attachment, we perform complexification-averaging and the
slow-flow averaging analysis. We elucidate the factors affecting the targeted
energy transfer from the main structure to the attachment, and use another time-scale to capture the evolution of the
system (super-slow flow) near the fixed point in the averaged phase space.
The optimal parameters are found by analyzing the system at super-slow time scale. We provide perturbation theoretic arguments along with numerical evidence for the validity of the model.
 In Section \ref{tmd}, we
describe a semi-analytical process to optimize a linear tuned mass damper
attached to one and two-degree-of-freedom system. We use
Lyapunov analysis to formulate the optimization problem, using energy
dissipated through the attachment as the metric. In Section \ref{results}, we
compare the vibration suppression performance of the two optimized
attachments, i.e. the NES and linear TMD, both attached to a
two-degree-of-freedom system. The results show that while NES performance is
sensitive to the energy of the impulse input, it is more robust than TMD in several scenarios.  In the Appendix, generalization of this framework for optimizing the
parameters of an essentially nonlinear attachment, coupled with an
n-degree-of-freedom system is provided.

\section{Dynamics of Multi-degree-of-freedom Primary System with NES}

\label{two-dof}
\begin{figure}[th]
\centering
\includegraphics[scale=0.5]{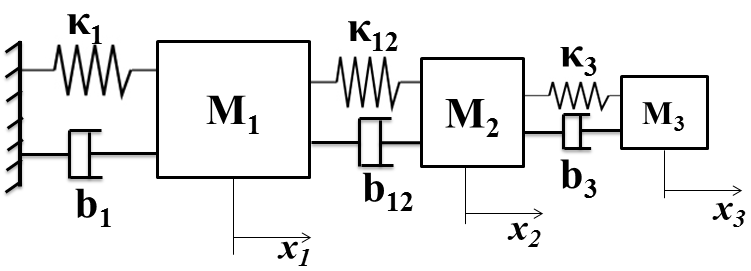} \caption{two-degree-of-freedom
system with NES.}%
\label{fig_2dof01}%
\end{figure}

A schematic of the two-degree-of-freedom system with NES is shown in
Fig. \ref{fig_2dof01}. Mass $M_{3}$ is assumed to be attached to mass $M_{2}$ using a
cubic spring with coefficient $\kappa_{3}$. The equations of motion for this
system are%

\begin{subequations}
\label{eq_2dof01}%
\begin{gather}
M_{1}\ddot{x}_{1}+b_{1}\dot{x}_{1}+b_{12}(\dot{x}_{1}-\dot{x}_{2})+\kappa
_{1}x_{1}+\kappa_{12}(x_{1}-x_{2})=0,\\
M_{2}\ddot{x}_{2}+b_{12}(\dot{x}_{2}-\dot{x}_{1})+b_{3}(\dot{x_{2}}-\dot
{x_{3}})+\kappa_{12}(x_{2}-x_{1})+\kappa_{3}(x_{2}-x_{3})^{3}=0,\\
M_{3}\ddot{x}_{3}+b_{3}(\dot{x}_{3}-\dot{x}_{2})+\kappa_{3}(x_{3}-x_{2}%
)^{3}=0.
\end{gather}
Defining the non-dimensional time%

\end{subequations}
\begin{equation}
\tau=\sqrt{\frac{\kappa_{1}}{M_{1}}}t, \label{eq_2dof02}%
\end{equation}
Equations (\ref{eq_2dof01} (a-c)) can now be written in a non-dimensional form
as
\begin{subequations}
\label{eq_2dof03}%
\begin{gather}
x_{1}^{\prime\prime}+2\zeta_{1}x_{1}^{\prime}+2\zeta_{12}(x_{1}^{\prime}%
-x_{2}^{\prime})+x_{1}+k_{12}(x_{1}-x_{2})=0,\\
\mu x_{2}^{\prime\prime}+2\zeta_{12}(x_{2}^{\prime}-x_{1}^{\prime})+2\zeta
_{3}(x_{2}^{\prime}-x_{3}^{\prime})+k_{12}(x_{2}-x_{1})+C(x_{2}-x_{3}%
)^{3}=0,\\
\epsilon x_{3}^{\prime\prime}+2\zeta_{3}(x_{3}^{\prime}-x_{2}^{\prime
})+C(x_{3}-x_{2})^{3}=0,
\end{gather}

where $^{\prime}$ denotes a derivative with respect to $\tau$ and $\mu
=\frac{M_{2}}{M_{1}}$, $\epsilon=\frac{M_{3}}{M_{1}}$, $\zeta_{1}=\frac{b_{1}%
}{2\sqrt{M_{1}\kappa_{1}}}$, $\zeta_{12}=\frac{b_{12}}{2\sqrt{M_{1}\kappa_{1}%
}}$ , $\zeta_{3}=\frac{b_{3}}{2\sqrt{M_{1}\kappa_{1}}}$ , $k_{12}=\frac
{\kappa_{12}}{\kappa_{1}}$ and $C=\frac{\kappa_{3}}{\kappa_{1}}$. 
The modal frequencies of the undamped primary system are $0.76$ rad/s and $2.63$ rad/s, while the corresponding mode shapes are $[0.65 \:\:\: 0.76]'$ and $[0.6 \:\:\: -0.8]'$, respectively.

\begin{table}[ptb]
\centering%
\begin{tabular}
[c]{cc}\hline
Parameter & Value\\\hline
$M_{1}$ [Kg] & 2200\\
$M_{2}$ [Kg] & 1400\\
$M_{3}$ [Kg] & 70\\
$\kappa_{1}$ [N/m] & $5.2\times10^{5}$\\
$\kappa_{12}$ [N/m] & $1.3\times10^{6}$\\
$\kappa_{3}$ [N/m] & $2.6\times10^{5}$\\
$b_{1}$ [Ns/m] & $5\times10^{2}$\\
$b_{12}$ [Ns/m] & $1\times10^{3}$\\
$b_{3}$ [Ns/m] & $50$%
\end{tabular}
\caption{Parameter values for the system in Fig. \ref{fig_2dof01} .}\label{tab_const01}%
\end{table}

The various modes of energy transfer energy transfer from single-degree-of-freedom primary system to NES were studied in Refs. \cite{vakphysD01, kerschen2005irreversible}. It has been shown that in the case of weakly damped primary systems, the energy transfer is mediated by resonance capture into, and escape from various periodic orbits of the \emph{undamped} system. 
This class of system has been shown to exhibit three main modes of energy transfer from the primary system to the NES, namely, targeted energy transfer through
\end{subequations}
\begin{enumerate}
\item fundamental transient resonance capture,

\item subharmonic transient resonance capture,
\begin{figure}[th]
\centering
\includegraphics[scale=0.5]{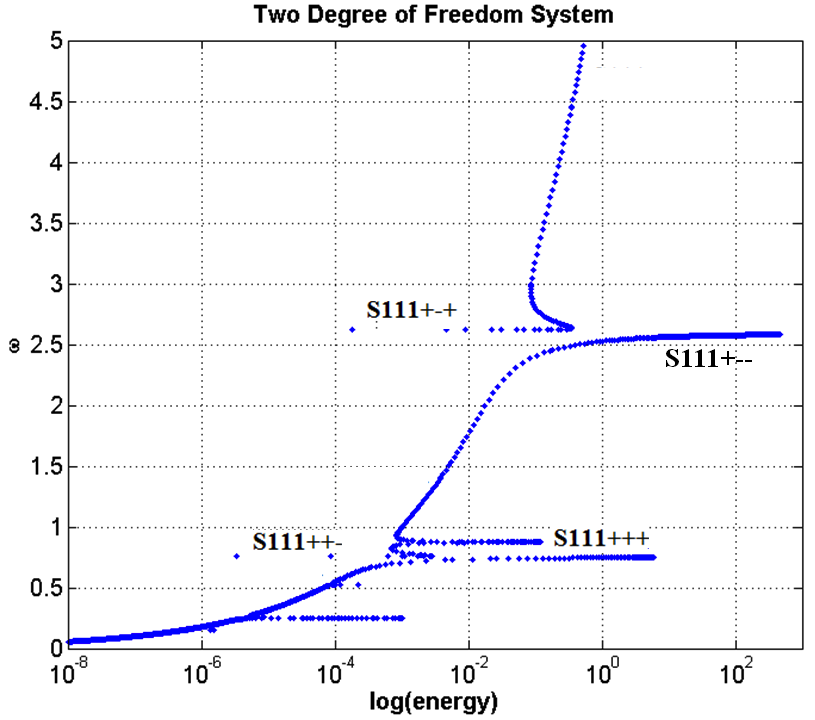}\caption{Frequency energy
plot showing periodic orbits of the system given in Eq. (\ref{eq_2dof03}).}%
\label{fig_2dof02}
\end{figure}

\item nonlinear beats.
\end{enumerate}
\begin{figure}[th!]
\centering
\includegraphics[width=.45\textwidth]{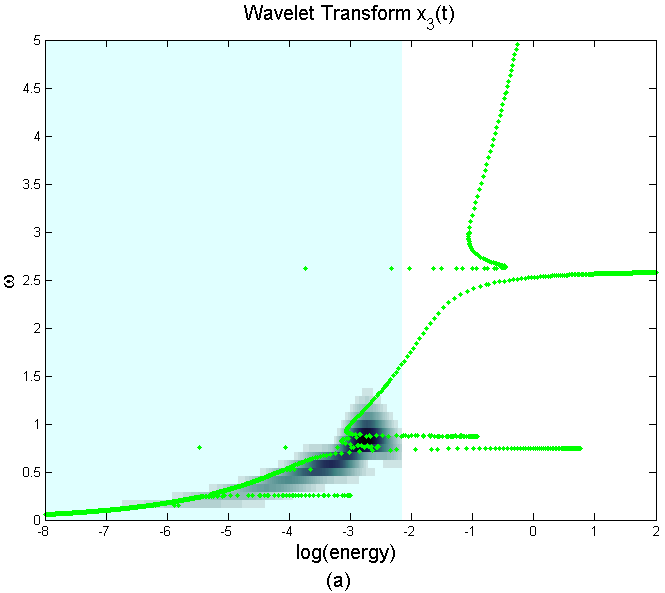}
\includegraphics[width=.45\textwidth]{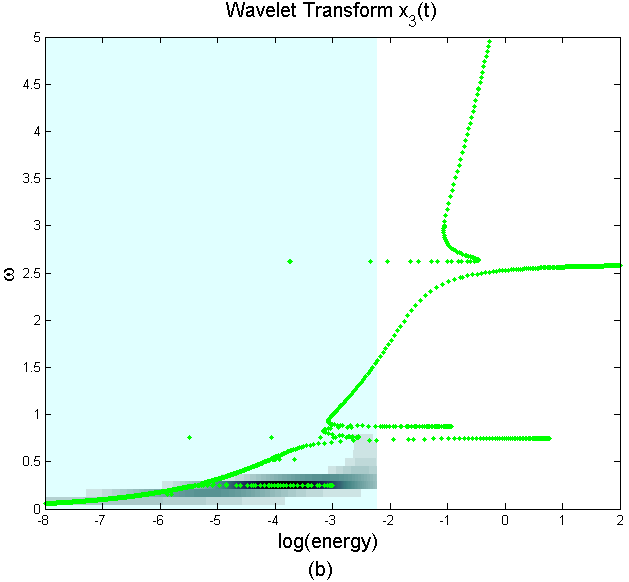}
\includegraphics[width=.45\textwidth]{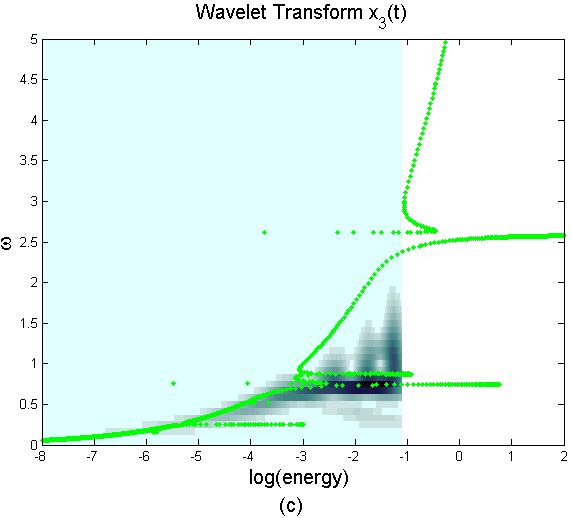} \caption{The three of modes of
targeted energy transfer graphically illustrated using wavelet transforms, (a)
Fundamental transient resonance capture, (b) Subharmonic transient resonance
capture, and (c) Nonlinear beats.}%
\label{fig_2dof03}%
\end{figure}

A comprehensive study of dynamics of a two-degree-of-freedom primary system with an attached NES was performed in Ref. \cite{kerschen2007theoretical}. Compared to the single-degree-of-freedom case, there are many more sequences of resonance transitions possible in this system, depending upon the initial condition. It is known \cite{vakbook} that the first two types of energy transfer modes cannot occur directly when the system is initially at rest. Therefore whenever the response of the system to an impulse is considered, the energy is always first transferred through nonlinear beats until one of the first two energy transfer modes is activated. 

Periodic orbits of the system of equations given in Eq. (\ref{eq_2dof03}) can
be found using non-smooth transformations \cite{pilipchuk02} along with the
numerical continuation software AUTO \cite{doedel1997auto}. Some of those orbits are shown in the frequency energy plot given in Fig.
\ref{fig_2dof02}.\newline

The fundamental transient resonance capture occurs when the system moves on the $S111+++$ or $S111+--$ branch as shown on the frequency-energy plot in Fig. \ref{fig_2dof02}. These branches have frequencies equal to the modal frequencies of the primary system. The symbols $+++$ indicate the phase of the three masses in the periodic orbit as plotted
in Fig. \ref{fig_2dof02}, i.e., $+++$ means that all three masses are in-phase
over the course of the system executing its periodic motion. Similarly $+--$
implies that the two primary masses are out-of-phase with each other. The numbers 111 in $S111$
denote the number of half sine waves the masses have in their half periods.
Thus in the $S111$ orbit, all three masses have a single half-sine wave in their
period. Fig. \ref{fig_2dof03} (a) uses wavelet transforms superimposed on
the frequency energy plot to illustrate targeted energy transfer mechanism while the system is on $S111+++$. Due to fundamental resonance capture, almost the entire energy of the system is localized in the NES. This makes it the most desirable mechanism for targeted energy transfer among the three mechanisms listed above  \cite{kerschen2007theoretical}.

Subharmonic energy transfer occurs when the initial energy given to the system is not enough to excite the fundamental transient resonance capture, and it is not as efficient as the fundamental resonant capture. One subharmonic energy transfer mode is illustrated in Fig. \ref{fig_2dof03} (b). 

Energy transfer through nonlinear beats occurs when the initial system energy is higher than the energy required to excite the fundamental energy transfer. In this mechanism first the NES undergoes a nonlinear beating phenomenon with the primary system and then gets attracted or \textquotedblleft captured\textquotedblright\ onto the fundamental energy transfer mode or one of the subharmonic modes.

When the system is given an impulse by disturbing one of the primary masses, (say) $M_1$, a multi-modal response is expected. Note that the out-of-phase $S111+--$ mode has a much higher energy threshold than the in-phase $S111+++$ mode. Hence, if one considers low to medium energy impulses, the optimal way to remove energy from the primary system is to get the system captured in the $S111+++$ mode, via the above-mentioned nonlinear beating phenomenon. In what follows, we analyze the dynamics of this scenario in detail, and obtain system parameters that result in near-optimal energy transfer for low to medium energy impulses.

\section{Optimization of targeted energy transfer}

\begin{figure}[th]
\centering
\includegraphics[width=0.75\textwidth]{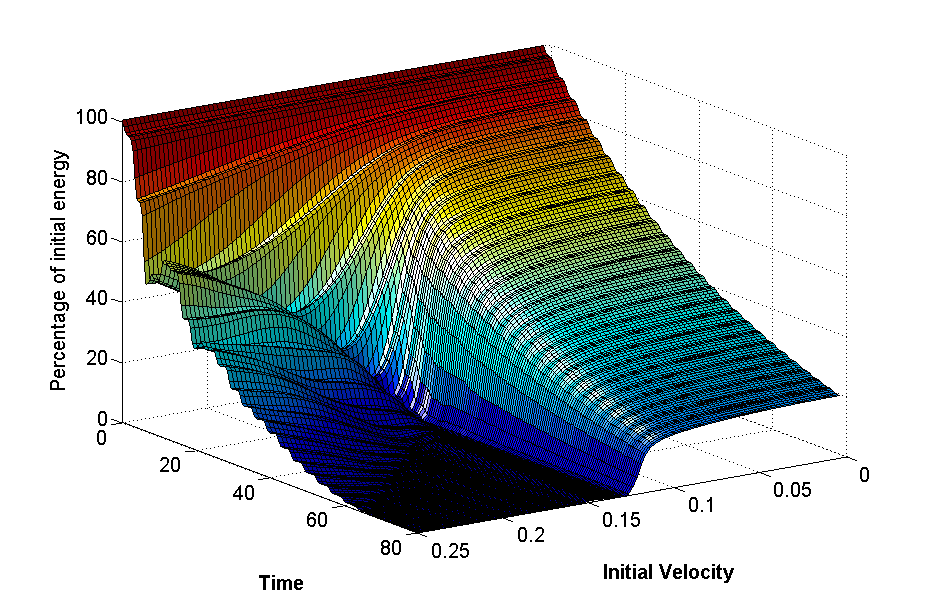} \caption{1 DOF system with NES: Percentage of initial
energy in the primary system with time for different starting velocities for a
one-degree-of-freedom system attached to the NES. The system is the one given
in \cite{sapsis_01}.}%
\label{fig_2dof04}%
\end{figure}

In a work by Sapsis et al. \cite{sapsis_01}, the energy transfer between a one-degree-of-freedom linear
oscillator coupled to a light nonlinear attachment was analyzed. It was shown
that optimal energy transfer from the primary linear oscillator to the
nonlinear attachment occurs when the initial energy of the system is close to
a particular homoclinic orbit of the underlying Hamiltonian system. Our analysis uses the strategy outlined in \cite{sapsis_01}, followed by dimensionality reduction to get explicit expressions for optimal parameters in the multiple-degree-of-freedom primary system.
\newline
For the system described by Eq.
(\ref{eq_2dof03}), the instantaneous energy stored in the primary system
(masses $M_{1}$ ans $M_{2}$) can be written as
\begin{equation}
\label{eq_2dof04}E_{inst}=\frac{1}{2}(x_{1}^{2}+x_{1}^{\prime2})+\frac{1}%
{2}\mu x_{2}^{\prime2}+\frac{1}{2}k_{12}(x_{1}-x_{2})^{2}.
\end{equation}
Consider the initial condition when only mass $M_{1}$ has a non-zero velocity
($v_{0}$). In that case the starting energy of the primary system is
\begin{equation}
\label{eq_2dof05}E_{start}=\frac{1}{2}v_{0}^{2}.
\end{equation}
The desired purpose of the NES is to remove energy from the primary system
irreversibly. One of the ways to quantify effectiveness of the NES used in
\cite{sapsis_01} is to look at the evolution of the quantity $\frac{E_{inst}%
}{E_{start}}$ for different starting velocities. This leads to the conclusion
that there is a threshold velocity above which the fundamental transient
resonance capture mode of TET is triggered and thus one can observe a dramatic
increase in the rate of drop of the ratio $\frac{E_{inst}}{E_{start}}$ due to
energy transfer from the primary system to the NES.
Fig. \ref{fig_2dof04} shows the change in instantaneous energy of the
primary system given in \cite{sapsis_01}. It can be observed that the rate of energy removal from the
primary system undergoes a jump around $v_{0}=0.11$. This is the starting
velocity which endows the system with the right amount of energy to be in the
neighborhood of the homoclinic orbit as discussed in \cite{sapsis_01}.\newline

\begin{figure}[th]
\centering
\includegraphics[width=0.75\textwidth]{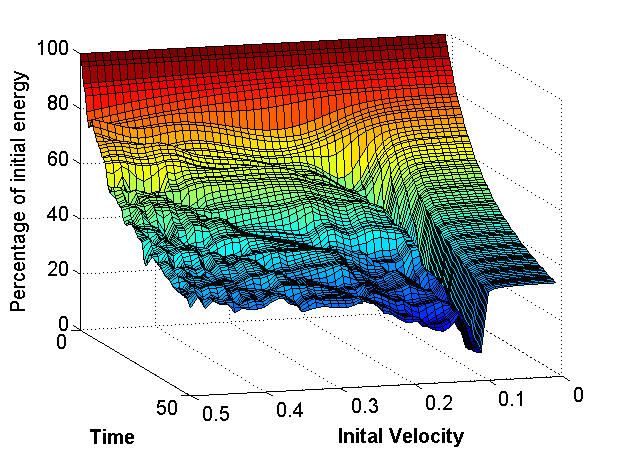} \caption{2 DOF system with NES: Percentage of initial
energy in the primary system with time for different starting velocities for
the system given in Fig. \ref{fig_2dof01}}%
\label{fig_2dof06}%
\end{figure}

For the system given by Eq. (\ref{eq_2dof03}), damping performance of the nonlinear attachment is shown in Fig. \ref{fig_2dof06}. The parameters of the system are given in Table \ref{tab_const01}. Unless otherwise noted, these parameters will be used throughout the paper. The rate of energy removal from the primary system has a sudden increase around $v_{0}=0.115$ and then the performance of the NES is sustained over a higher initial velocity range. Hence for right choice of initial velocity, the nonlinear attachment still acts like an energy sink for the two-degree-of-freedom primary system.

\subsection{Study Using Complexification-Averaging}

For approximate analytical study of the system given in Eq. (\ref{eq_2dof03}),
the technique of complexification-averaging used in \cite{vakphysD01,
sapsis_01} is employed. We study the dynamics in the vicinity of the fundamental transient
resonance $S111+++$ , hence the fast frequency for all three masses is $\omega$. As mentioned erlier, we only consider low to medium energy initial impulses, i.e. case where the in-phase fundamental frequency $\omega$ is not significantly dependent on the initial energy. 
Introducing the new variables%

\begin{subequations}
\begin{gather}
\psi_{1}=x_{1}^{\prime}(\tau)+\omega jx_{1}(\tau),\\
\psi_{2}=x_{2}^{\prime}(\tau)+\omega jx_{2}(\tau),\\
\psi_{3}=x_{3}^{\prime}(\tau)+\omega jx_{3}(\tau),
\end{gather}
and substituting%

\end{subequations}
\begin{equation}
\psi_{i}=\phi_{i}e^{j\omega\tau},\;i=1,2,3,
\end{equation}
Equation (\ref{eq_2dof03}) can be averaged over the fast time scale $\tau$. This yields%

\begin{subequations}
\label{eq_2dof08}%
\begin{gather}
\phi_{1}^{\prime}+\left(  \zeta_{1}+\frac{j\omega}{2}-\frac{j}{2\omega
}\right)  \phi_{1}+\zeta_{12}(\phi_{1}-\phi_{2})-\frac{jk_{12}}{2\omega}%
(\phi_{1}-\phi_{2})=0,\\
\mu\phi_{2}^{\prime}+\zeta_{12}(\phi_{2}-\phi_{1})+\zeta_{3}(\phi_{2}-\phi
_{3})+\frac{j\omega}{2}\mu\phi_{2}-\frac{jk_{12}}{2\omega}(\phi_{2}-\phi
_{1})-\frac{3jC}{8\omega^{3}}|\phi_{2}-\phi_{3}|^{2}(\phi_{2}-\phi_{3})=0,\\
\epsilon\phi_{3}^{\prime}+\zeta_{3}(\phi_{3}-\phi_{2})+\frac{j\omega}%
{2}\epsilon\phi_{3}-\frac{3jC}{8\omega^{3}}|\phi_{3}-\phi_{2}|^{2}(\phi
_{3}-\phi_{2})=0.
\end{gather}
Now we introduce the new variables%

\end{subequations}
\begin{align}
u_{1}=\phi_{1}-\phi_{2},\nonumber\\
u_{2}=\phi_{2}-\phi_{3},\nonumber\\
u_{3}=\phi_{1}+\mu\phi_{2}+\epsilon\phi_{3},
\end{align}

where $u_{1}$ represents the relative displacement between the first mass and
the second mass, $u_{2}$ the relative displacement between the second mass and
the NES and $u_{3}$ represents the motion of the center of mass of the system.
The equations for $u_{i}$ can be derived from Eq. (\ref{eq_2dof08}) and using
the relations
\begin{subequations}
\begin{gather}
\phi_{1}=\frac{u_{3}+u_{1}(\mu+\epsilon)+u_{2}\epsilon}{1+\mu+\epsilon},\\
\phi_{2}=\frac{u_{3}-u_{1}+u_{2}\epsilon}{1+\mu+\epsilon},\\
\phi_{3}=\frac{u_{3}-u_{1}-u_{2}(1+\mu)}{1+\mu+\epsilon}.
\end{gather}
The equations are%

\end{subequations}
\begin{subequations}
\begin{gather}
u_{1}^{\prime}+c_{11}u_{1}+c_{12}u_{2}+c_{13}u_{3}+\frac{3jC}{8\mu\omega^{3}%
}|u_{2}|^{2}u_{2}=0,\\
u_{2}^{\prime}+c_{21}u_{1}+c_{22}u_{2}-\frac{3jC(\mu+\epsilon)}{8\mu
\epsilon\omega^{3}}|u_{2}|^{2}u_{2}=0,\\
u_{3}^{\prime}+c_{31}u_{1}+c_{32}u_{2}+c_{33}u_{3}=0,
\end{gather}\label{eq_2dof15}
where%
\end{subequations}

\begin{gather}
c_{11}=-\frac{j(\mu^{2}+k_{12}(1+\mu)(1+\epsilon+\mu)+2j\mu^{2}\zeta_{1}%
\omega+2j\zeta_{12}\omega+4j\mu\zeta_{12}\omega+2j\mu^{2}\zeta_{12}\omega
-\mu\omega^{2}-\mu^{2}\omega^{2}+\epsilon N_{1})}{2\mu(1+\epsilon+\mu)\omega
},\label{eq_2dof16}\\
N_{1}=\mu+2j\mu\zeta_{1}\omega+2j\zeta_{12}\omega+2j\mu\zeta12\omega-\mu
\omega^{2},\\
c_{12}=\frac{-2(1+\mu)\zeta_{3}\omega+\epsilon(-2\zeta_{3}\omega+\mu
(-j+2\zeta_{1}\omega))}{2\mu(1+\epsilon+\mu)\omega},\;\;c_{13}=\frac
{-j+2\zeta_{1}\omega}{2(1+\epsilon+\mu)\omega},\\
c_{21}=\frac{j(k_{12}+2j\zeta_{12}\omega)}{2\mu\omega},\;\;c_{22}=\frac
{\zeta_{3}}{\epsilon}+\frac{\zeta_{3}}{\mu}+\frac{j\omega}{2},\\
c_{31}=\frac{(\epsilon+\mu)(-j+2\zeta_{1}\omega)}{2(1+\epsilon+\mu)\omega
},\;\;c_{32}=\frac{\epsilon(-j+2\zeta_{1}\omega)}{2(1+\epsilon+\mu)\omega
},\;\;c_{33}=\frac{2\zeta_{1}\omega+j(-1+(1+\epsilon+\mu)\omega^{2}%
)}{2(1+\epsilon+\mu)\omega},
\end{gather}

The initial conditions for $u_{i}$ are%
\begin{subequations}
\begin{gather}
u_{1}(0)=\phi_{1}(0)-\phi_{2}(0)=v_{0},\label{eq_2dof17}\\
u_{2}(0)=\phi_{2}(0)-\phi_{3}(0)=0,\\
u_{3}(0)=\phi_{1}(0)+\mu\phi_{2}(0)+\epsilon\phi_{3}(0)=v_{0}.
\end{gather}\label{eq_2dof17}

The instantaneous energy in the primary system in terms of the `slow'
variables can be written as%

\begin{gather}
E_{p}=\frac{1}{2}(x_{1}^{2}+x_{1}^{\prime2})+\frac{1}{2}\mu x_{2}^{\prime
2}+\frac{1}{2}k_{12}(x_{1}-x_{2})^{2}\approx\nonumber\\
\frac{1}{2}\left(  \operatorname{Re}[\phi_{1}e^{j\omega\tau}]^{2}+\frac
{1}{\omega^{2}}\operatorname{Im}[\phi_{1}e^{j\omega\tau}]^{2}\right)
+\frac{1}{2}\mu(\operatorname{Re}[\phi_{2}e^{j\omega\tau}])^{2}+\frac
{1}{2\omega^{2}}k_{12}(\operatorname{Im}[u_{1}e^{j\omega\tau}])^{2}.
\label{eq_2dof10}%
\end{gather}
The instantaneous energy stored in the NES can be written as%

\begin{gather}
E_{NES}=\frac{1}{2}\epsilon x_{3}^{\prime2}+\frac{1}{4}C(x_{2}-x_{3}%
)^{4}\approx\nonumber\\
\frac{1}{2}\epsilon({\operatorname{Re}}[\phi_{3}e^{j\omega\tau}])^{2}+\frac
{1}{4\omega^{4}}C(\operatorname{Im}[u_{2}e^{j\omega\tau}])^{4}.
\end{gather}
\label{eq_2dof11}

The energy dissipated by the NES can be written as
\begin{gather}
E_{DISS}=\int_{0}^{\tau}{\zeta_{3}(x_{2}^{\prime}-x_{3}^{\prime})^{2}}%
d\tau\approx\int_{0}^{\tau}{\zeta_{3}(\operatorname{Re}[u_{2}e^{j\omega\tau
}])^{2}}d\tau=\nonumber\\
\int_{0}^{\tau}{\zeta_{3}(Re[u_{2}]^{2}\cos^{2}{\omega\tau}%
+\operatorname{Im}[u_{2}]^{2}\sin^{2}{\omega\tau}-\operatorname{Re}%
[u_{2}]\operatorname{Im}[u_{2}]\sin{2\omega\tau})}d\tau, \label{eq_2dof12}%
\end{gather}
and omitting frequencies faster than $\omega$ as argued in \cite{sapsis_01},
one can approximately write
\begin{equation}
E_{DISS}\approx\int_{0}^{\tau}{\zeta_{3}(|u_{2}(\tau)|)^{2}}d\tau.
\label{eq_2dof13}%
\end{equation}

As noted in \cite{sapsis_01}, the performance of the NES is closely tied to
the variable on which $E_{DISS}$ depends, which in the system under
consideration happens to be $u_{2}$. Adapting the principle of optimal dissipation from
\cite{sapsis_01}, we claim that Enhanced TET in the system is realized
when $|u_{2}(\tau)|$ exhibits large amplitudes, especially during the initial
phase of motion when energy is highest.

\begin{figure}[th]
\centering
\includegraphics[scale=0.4]{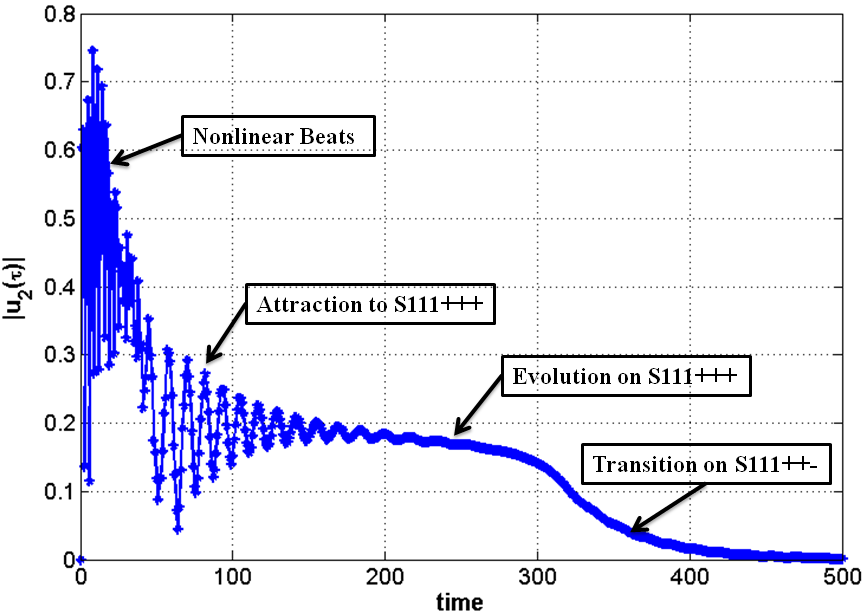}\caption{Evolution of $|u_{2}(\tau)|$
with different energy transfer regimes for a velocity input $v_{0} = 0.3$}\label{fig_2dof07}
\end{figure}
With Eqs. (\ref{eq_2dof15}) and the initial conditions (\ref{eq_2dof17}), the
slow time variables $u_{i}$ can be integrated forward in time and different
energy transfer regimes between the primary system and NES for different
initial velocities can be identified. Fig. \ref{fig_2dof07} shows a typical curve of
$|u_{2}(\tau)|$ for a particular starting non-zero velocity for the mass
$M_{1}$, all other initial conditions being zero. As mentioned earlier, since
fundamental transient resonance capture and subharmonic orbits cannot get
excited directly by impulsive inputs, there is a regime of nonlinear beats in
the initial stages of the response. Recall from Fig. \ref{fig_2dof06}, that
$v_{0}=0.115$ is close to the transition point where the NES performance
improves. The evolution for the energies in the primary system and the NES for
three different initial velocities is given in Fig.
\ref{fig_2dof08}. For $v_{0}=.05$, there is no energy transfer to the NES. For $v_{0}=.2$, some energy is transferred irreversibly to the NES. However, the most efficient energy transfer among the three cases is for $v_{0}=.115$, where almost all energy from the primary system is transferred to the NES within the given time. We will see shortly that the underlying mechanism for this behavior is the existence of a homoclinic orbit in the undamped system at super-slow time scale. If the system is given less energy than the energy of this orbit, efficient energy transfer does not occur, as is seen for $v_{0}=.05$. If the initial energy of the system is more than the energy of this homoclinic orbit, the system gets captured via nonlinear beats into the $S111+++$ orbit and eventually transfers all the energy to NES. If the initial energy is close to the energy of homoclinic orbit, i.e. the 'optimal' case, all energy is transferred to the NES in a single cycle of the nonlinear beats, i.e. this is a degenerate case at the boundary of the first two cases.

\begin{figure}[th]
\centering
\includegraphics[scale=0.5]{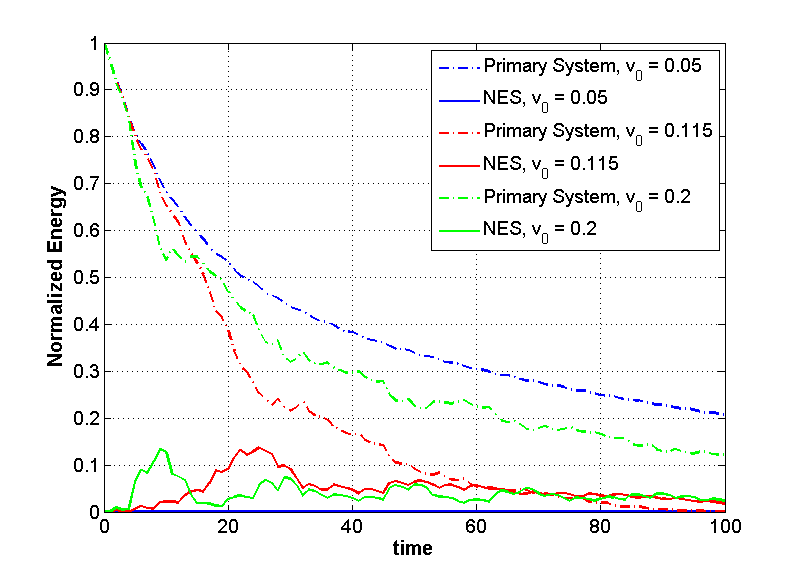}\caption{Evolution of
normalized energy in the primary system and NES for different starting
velocities.}\label{fig_2dof08}
\end{figure}

\subsection{Analysis of Super-Slow Flow Equations}

The solution for Eq. (\ref{eq_2dof15}c)  is%
\end{subequations}
\begin{equation}
u_{3}=u_{3}(0)e^{-c_{33}\tau}-c_{31}\int_{0}^{\tau}{u_{1}e^{c_{33}(t_{1}%
-\tau)}dt_{1}}-c_{32}\int_{0}^{\tau}{u_{2}e^{c_{33}(t_{1}-\tau)}dt_{1}.}
\label{eq_2dof18}%
\end{equation}
Now for parameters under consideration (for low to medium energy impulses), we find
\begin{equation}
(1+\mu+\epsilon)\omega^{2}\approx1, \label{eq_2dof19}%
\end{equation}
one can write
\begin{equation}
c_{33}=\frac{\zeta_{1}}{1+\mu+\epsilon}. \label{eq_2dof20}%
\end{equation}
As the damping coefficients $\zeta_{1}$ and $\zeta_{3}$ are small, one can assume%

\begin{subequations}
\begin{align}
\zeta_{1}  &  =\epsilon\widehat{\zeta}_{1},\\
\zeta_{3}  &  =\epsilon\widehat{\zeta}_{3}.
\end{align}

Further assuming that $t_{1}-\tau=O(1/\epsilon^{1/2})$, Eq. \ref{eq_2dof18}
can be re-written as%

\end{subequations}
\begin{equation}
u_{3}\approx u_{3}(0)e^{-c_{33}\tau}-c_{31}\overline{u_{1}}\tau e^{-c_{33}%
(\overline{t}-\tau)}+O(\epsilon), \label{eq_2dof22}%
\end{equation}
where $\overline{u}_{1}$ is the value of $u_{1}$ at $t=\overline{t}$ in the
interval $[0,\tau]$ according to the mean-value theorem. For simplicity, if we assume $\overline{t}\approx\tau$, then Eq. (\ref{eq_2dof22}) can be further
simplified as
\begin{equation}
u_{3}\approx u_{3}(0)e^{-c_{33}\tau}-c_{31}u_{1}(\tau)\tau+O(\epsilon).
\label{eq_2dof22b}%
\end{equation}
For obtaining the super-slow flow equations, a new scaling is introduced:
\begin{subequations}
\begin{align}
u_{1}  &  =\epsilon^{1/2}z_{1},\\
u_{2}  &  =\epsilon^{1/2}z_{2},\\
u_{3}  &  =\epsilon^{1/2}z_{3}.
\end{align}
Using this scaling and Eq. (\ref{eq_2dof22b}), Eqs. (\ref{eq_2dof15}(a),\ref{eq_2dof15}(b)) can be re-written as%

\end{subequations}
\begin{subequations}
\begin{gather}
z_{1}^{\prime}+\widehat{c}_{11}z_{1}+\widehat{c}_{13}\left(  z_{3}%
(0)e^{-c_{33}\tau}-c_{31}z_{1}\tau\right)  +O(\epsilon)=0,\\
z_{2}^{\prime}+\widehat{c}_{21}z_{1}+\widehat{c}_{22}z_{2}-\frac{3jC}%
{8\omega^{3}}|z_{2}|^{2}z_{2}+O(\epsilon)=0,
\end{gather}\label{eq_2dof24}
where%

\end{subequations}
\begin{gather}
\widehat{c}_{11}=\frac{-j\mu}{2\omega(1+\mu+\epsilon)}+\frac{j\omega}%
{2}+\left(  \zeta_{12}-\frac{jk_{12}}{2\omega}\right)  \left(  1+\frac{1}{\mu
}\right)  ,\label{eq_2dof24_b}\\
\widehat{c}_{13}=\frac{-j}{2(1+\epsilon+\mu)\omega},\\
\widehat{c}_{21}=\frac{j(k_{12}+2j\zeta_{12}\omega)}{2\mu\omega}%
,\;\;\widehat{c}_{22}=\widehat{\zeta}_{3}+\frac{j\omega}{2}.
\end{gather}
\begin{figure}[th!]
\centering
\includegraphics[scale=0.5]{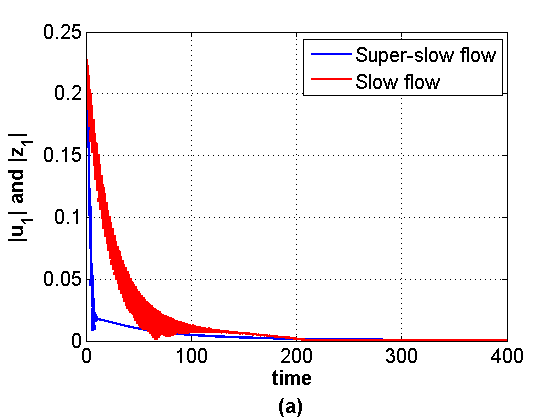} \includegraphics[scale=0.5]{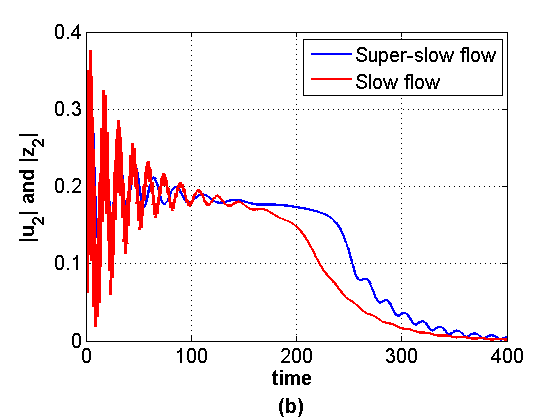}
\caption{Comparison of the solutions of the slow and super-slow equations.}%
\label{fig_2dof09}%
\end{figure}

\begin{figure}[th]
\centering
\includegraphics[scale=0.5]{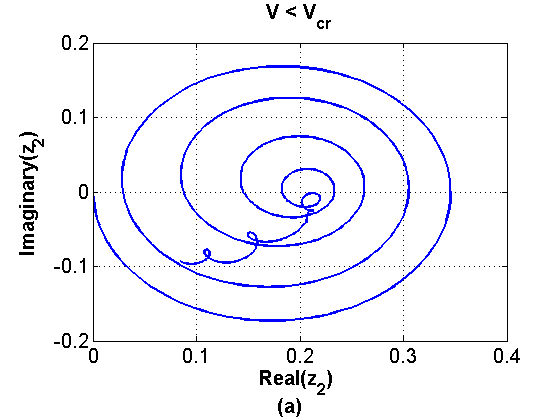} \includegraphics[scale=0.5]{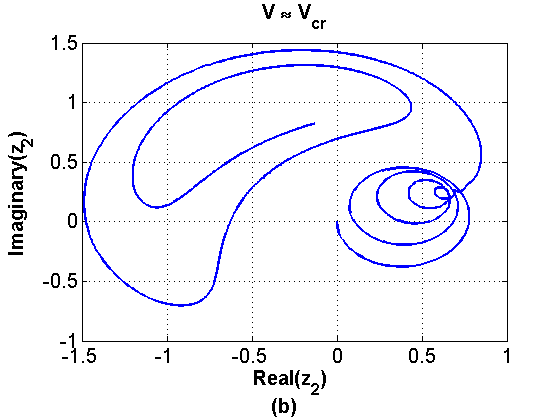}
\caption{Phase portrait of $z_{2}$ for two different initial velocities using the super-slow flow equations.}%
\label{fig_2dof10}%
\end{figure}

\begin{figure}[th]
\centering
\includegraphics[scale=0.45]{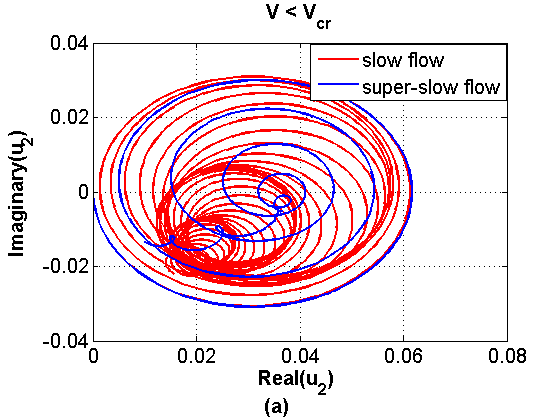} \includegraphics[scale=0.45]{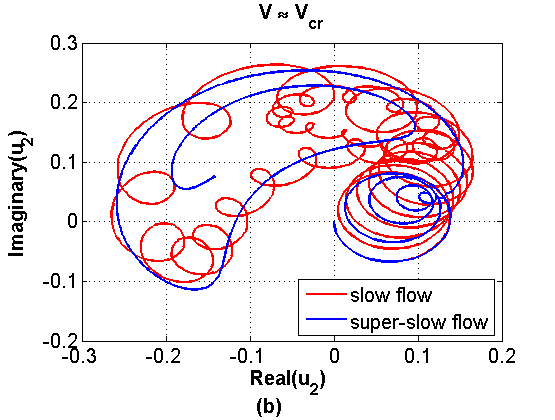}
\caption{Comparison of the phase portrait of $z_{2}$ for two different initial
velocities obtained using slow and super-slow flow equations.}%
\label{fig_2dof11}%
\end{figure}

Equations (\ref{eq_2dof24}) (a) and (b) are the super-slow flow equations. The
comparison of numerical integration for Eq. (\ref{eq_2dof24}) (a) and (b)
(super-slow) as well as for Eq. (\ref{eq_2dof15}) (a) and (b) (slow) are given
in Fig. \ref{fig_2dof09}. The $z_2$ phase-portrait of super-slow flow equations for two different initial velocities is shown in Fig. \ref{fig_2dof10}. We discuss these two cases later in the section. In Fig. \ref{fig_2dof11}, the phase-portrait solution of slow-flow equations is superimposed on that of the super-slow flow response. The super-slow-flow trajectories can be seen as low-pass filtered versions of slow-flow trajectories.

The solution for Eq. (\ref{eq_2dof24}(a)) (up to $O(\epsilon
^{0})$) can be written as%
\begin{equation}
z_{1}=z_{1c}+z_{1p}, \label{eq_2dof25}%
\end{equation}
where $z_{1c}$ is the complimentary solution and $z_{1p}$ is the particular
solution. For the complimentary solution only the homogenous part of the Eq.
(\ref{eq_2dof24}(a)) is considered, and thus the solution can be written as%

\begin{equation}
z_{1c}=z_{c0}e^{-\widehat{c}_{11}\tau+\widehat{c}_{13}c_{31}\tau^{2}/2},
\label{eq_2dof25a}%
\end{equation}
where $z_{c0}$ is the constant of integration. Looking at the form of the
non-homogenous term in Eq. (\ref{eq_2dof24}(a)) , the particular solution is
assumed to be
\begin{equation}
z_{1p}=Ae^{-c_{33}\tau}. \label{eq_2dof25b}%
\end{equation}
Substituting Eq. (\ref{eq_2dof25b}) into Eq. (\ref{eq_2dof24}(a))  and
simplifying, we obtain
\begin{equation}
A=\frac{-\widehat{c}_{13}z_{3}(0)}{\widehat{c}_{11}-\widehat{c}_{13}c_{31}%
\tau-c_{33}}. \label{eq_2dof25c}%
\end{equation}
Recall that, $A$ was assumed to be constant, i.e., independent of $\tau$ in
Eq. (\ref{eq_2dof25b}) and yet Eq. (\ref{eq_2dof25c}) gives a time dependence for
$A$ and thus there is a contradiction. However, the particular solution decays
exponentially (even for constant $A$). The $\tau$ dependence for $A$ is in
the denominator, and this dependence is dominated by the exponential decay, hence we  
neglect its derivative. In that case, the solution for Eq.
(\ref{eq_2dof24}(a)) up to $O(\epsilon^{0})$ can be approximated by%

\begin{equation}
z_{1}=z_{c0}e^{-\widehat{c}_{11}\tau+\widehat{c}_{13}c_{31}\tau^{2}/2}%
-\frac{\widehat{c}_{13}z_{3}(0)}{\widehat{c}_{11}-\widehat{c}_{13}c_{31}%
\tau-c_{33}}e^{-c_{33}\tau}, \label{eq_2dof25d}%
\end{equation}
where
\begin{equation}
z_{c0}=z_{1}(0)-A. \label{eq_2dof25e}%
\end{equation}

\begin{figure}[th]
\centering
\includegraphics[scale=0.5]{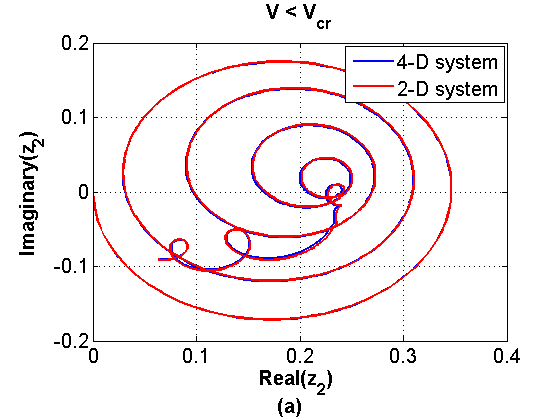}
\includegraphics[scale=0.5]{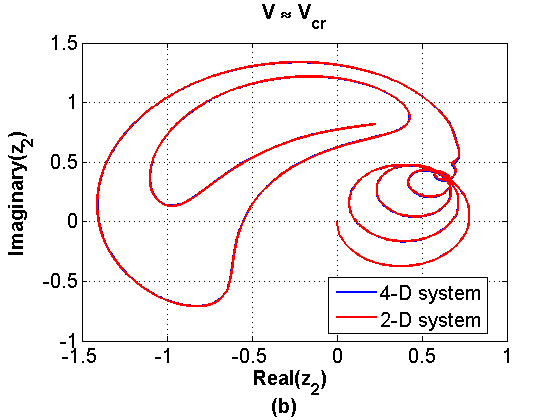} \caption{Comparison of the phase
portrait of $z_{2}$ for two different initial velocities with the 4 and 2
dimensional super-slow systems.}%
\label{fig_2dof12}%
\end{figure}

Now we substitute the solution for $z_{1}$ from Eq. (\ref{eq_2dof25d}) into Eq. (\ref{eq_2dof24}(b)) to obtain a single differential equation in $z_2$.  We refer to
this equation as the 2-D (two-dimensional) system. We refer to the coupled Eqs. (\ref{eq_2dof24}(a)-(b)) as the 4-D (four-dimensional) system.  These two systems can be separately integrated in time to
allow a comparison of their responses, and assessment of validity of the dimensionality reduction from 4-D to 2-D. In Fig. \ref{fig_2dof12}, we plot the responses of
both of these systems. It is clear from this figure that the response on the super-slow timescale is well reproduced by the 2-D system, for both cases of initial velocity shown in Fig. \ref{fig_2dof12}. Hence, using a series of averaging operations and approximations, we have reduced the system to one degree of freedom on the super-slow time-scale. 

In Fig. \ref{fig_2dof12}(a), the system is initialized with lower initial velocity than the critical velocity, denoted by $v_{cr}$. The trajectory starts at origin (since $z_2=0$ when initial impulse is given to $M_1$), grows in size for some time and then decays. The trajectory shown in Fig.  \ref{fig_2dof12}(b), initiated with close to critical velocity, makes large transitions in the phase space after an initial nonlinear beating phase, and gets captured in $S111+++$ resonance. 

We claim that for critical initial velocity $v_{cr}$, the undamped system has an homoclinic orbit, and this leads to almost complete energy transfer in one nonlinear beating cycle for the damped system. Hence, this critical case acts as the boundary between the low-energy case where there is no energy transfer to NES, and the high-energy case where the energy transfer is less than optimal. To support this claim, we study the undamped dynamics of the super-slow system in more detail. The sole equation for $z_{2}$ can be written as%

\begin{equation}
z_{2}^{\prime}+\left(  \widehat{\zeta}_{3}+\frac{j\omega}{2}\right)
z_{2}+\left(  \frac{jk_{12}-2\zeta_{12}\omega}{2\mu\omega}\right)
z_{co}e^{-\widehat{c}_{11}\tau+\widehat{c}_{13}c_{31}\tau^{2}/2}+\left(
\frac{jk_{12}-2\zeta_{12}\omega}{2\mu\omega}\right)  Ae^{-c_{33}\tau}%
-\frac{3jC}{8\omega^{3}}|z_{2}|^{2}z_{2}+O(\epsilon)=0. \label{eq_2dof26}%
\end{equation}
To find a conservative equivalent of Eq.(\ref{eq_2dof26}), the damping terms
$\widehat{\zeta}_{3}$, $\zeta_{12}$ and $\zeta_{1}$ are assumed to be small
and can be neglected. Considering then terms up to $O(\epsilon^{0})$ Eq.
(\ref{eq_2dof26}) reduces to
\begin{equation}
z_{2}^{\prime}+\frac{j\omega}{2}z_{2}-\frac{3jC}{8\omega^{3}}|z_{2}|^{2}%
z_{2}=z_{3}(0)(C_{1}+jC_{2}), \label{eq_2dof27}%
\end{equation}
where
\begin{gather}
C_{1}=-\frac{k_{12}}{2\mu\omega}\left(  e^{-\alpha\tau^{2}}(-Ap_{2}%
sin(\beta\tau))-A_{2}\right)  ,\label{eq_2dof27_a}\\
C_{2}=-\frac{k_{12}}{2\mu\omega}\left(  e^{-\alpha\tau^{2}}(Ap_{2}%
cos(\beta\tau))+A_{1}\right)  ,\\
\alpha=\frac{\mu}{8\omega^{2}(1+\mu+\epsilon)^{2}}\;,\;\;\beta=\frac{\mu
}{2\omega(1+\mu+\epsilon)}-\frac{\omega}{2}+\frac{k_{12}}{2\omega}\left(
1+\frac{1}{\mu}\right)  ,\\
A_{1}=\frac{b}{a^{2}+b^{2}}\;,\;\;A_{2}=\frac{a}{a^{2}+b^{2}},\;\;,\;\;Ap_{2}%
=1-\frac{\mu}{-\mu^{2}-k_{12}(1+\mu)(1+\mu+\epsilon)+\omega^{2}\mu
(1+\mu+\epsilon)},\\
a=\frac{\mu\tau}{2\omega(1+\mu+\epsilon)},\;\;b=-\mu-k_{12}\frac{1+\mu}{\mu
}(1+\mu+\epsilon)+\omega^{2}(1+\mu+\epsilon).
\end{gather}
Now, assuming
\begin{equation}
z_{2}=Ne^{j\delta}, \label{eq_2dof28}%
\end{equation}
and separating the real and imaginary parts of Eq. (\ref{eq_2dof27}) the
following equations for the magnitude $N$ and phase $\delta$ can be obtained,%

\begin{subequations}
\begin{gather}
N^{\prime}=z_{3}(0)(C_{1}cos(\delta)+C_{2}sin(\delta)),\\
\delta^{\prime}+\frac{\omega}{2}-\frac{3C}{8\omega^{3}}N^{2}=\frac{z_{3}%
(0)}{N}(C_{2}sin(\delta)-C_{1}cos(\delta)).
\end{gather}
Alternatively, if one assumes, $z_{2}=x+jy\;$, then the equations for $x$ and
$y$ are given by,%

\end{subequations}
\begin{subequations}
\begin{gather}
x^{\prime}=\frac{\omega y}{2}-\frac{3C}{8\omega^{3}}(x^{2}+y^{2}%
)y+z_{3}(0)C_{1},\\
y^{\prime}=-\frac{\omega x}{2}+\frac{3C}{8\omega^{3}}(x^{2}+y^{2}%
)y+z_{3}(0)C_{2}.
\end{gather}
Both $C_{1}$ and $C_{2}$ are functions of time. For further simplification, if
they are replaced by their mean values over the $\tau$ range of interest, then%

\end{subequations}
\begin{subequations}
\begin{gather}
x^{\prime}=\frac{\omega y}{2}-\frac{3C}{8\omega^{3}}(x^{2}+y^{2}%
)y+z_{3}(0)\widehat{C}_{1},\\
y^{\prime}=-\frac{\omega x}{2}+\frac{3C}{8\omega^{3}}(x^{2}+y^{2}%
)x+z_{3}(0)\widehat{C}_{2}.
\end{gather}\label{eq_2do31}

Observing Eq. (\ref{eq_2do31}), the following Hamiltonian can be written%

\end{subequations}
\begin{equation}
h=\frac{\omega}{4}|z_{2}|^{2}-\frac{3C}{8\omega^{3}}\frac{|z_{2}|^{4}}%
{4}-z_{3}(0)\widehat{C}_{2}\frac{z_{2}+z_{2}^{\ast}}{2}-z_{3}(0)j\widehat
{C}_{1}\frac{z_{2}-z_{2}^{\ast}}{2}, \label{eq_2dof32}%
\end{equation}
where the asterisk ($^{\ast}$) denotes complex conjugate. Since the Hamiltonian is constant of motion, Eqs. \ref{eq_2do31} (a) and (b) can
be reduced to a single equation,
\begin{gather}
a^{\prime}=2\sqrt{f(a,z_{3}(0))},\label{eq_2dof33}\\
a=N^{2},\qquad f(a,z_{3}(0))=z_{3}(0)^{2}a-\frac{1}{\widehat{C}_{1}%
^{2}+\widehat{C}_{2}^{2}}\left(  \frac{\omega}{4}a-\frac{3Ca^{2}}{32\omega
^{3}}\right)  ^{2}. %
\end{gather}

\begin{figure}[th]
\centering
\includegraphics[width=.6\textwidth]{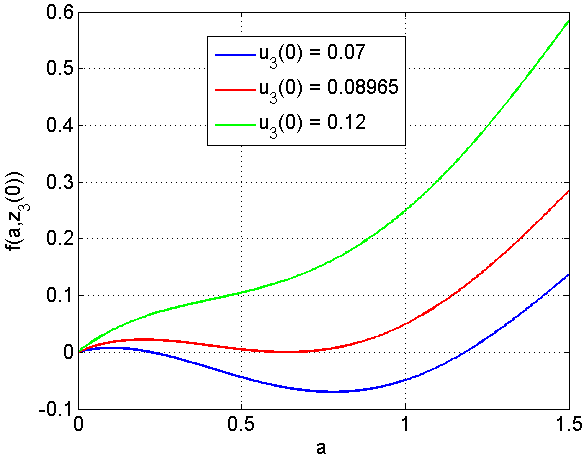} \caption{Roots of the equation
$f(a,z_{3}(0)) = 0\;$ .}%
\label{fig_2dof14}%
\end{figure}

Note that Eq. (\ref{eq_2dof33}) is a peculiar one, since it has the initial conditions on the right hand side. Hence, while $z_2(0)=0$ is needed for the case of interest (corresponding to the case where initial impulse is given to $M_1$), $z_3(0)$ acts as a tunable parameter in one-degree-of-freedom system given by Eq. (\ref{eq_2do31}) or Eq. (\ref{eq_2dof33}).
 In Fig.  \ref{fig_2dof14}, plot for the function $f(a,z_{3}(0))$ for
different values of $u_{3}(0) = \epsilon^{\frac{1}{2}} z_{3}(0)$ is shown. The system given by Eq. (\ref{eq_2dof33}) has four real fixed points in general. For a specific value of $z_3(0)$ (`the critical case'), two real roots of the function $f(a,z_{3}(0))$ coincide, similar to the observation made in Ref. \cite{sapsis_01}. This value can be found out by setting the derivative of $f(a,z_3(0))$ with respect to $a$ equal to zero, i.e. ,
\begin{equation}
f^{\prime}(a,z_3(0)) = z_3(0)_{cr}^2  - \frac{2}{\widehat{C}_1^2 + \widehat{C}_2^2} \left( \frac{\omega}{4} - \frac{6Ca}{32 \omega^3}\right)\left( \frac{\omega a}{4} - \frac{6Ca^2}{32 \omega^3}\right) = 0\;,
\label{eq_2dof34}
\end{equation}
where $z_3(0)_{cr}$ is the critical value. Also, since the fixed point is a root of the equation $f(a,z_3(0)) = 0$, one can write,

\begin{equation}
z_3(0)_{cr}^2 = \frac{a}{\widehat{C}_1^2 + \widehat{C}_2^2}\left( \frac{\omega}{4} - \frac{6Ca}{32 \omega^3}\right)^2 \;.
\label{eq_2dof35}
\end{equation}
Using Eq. (\ref{eq_2dof34}) and (\ref{eq_2dof35}), the expression for $z_3(0)_{cr}$ can be written as,

\begin{equation}
z_3(0)_{cr} = \frac{\omega^3}{9}\sqrt{\frac{2}{C (\widehat{C}_1^2 + \widehat{C}_2^2)}} \;.
\label{eq_2dof36}
\end{equation}

Therefore, the expression for the critical initial velocity given to $M_1$ is,
\begin{equation}
v_{cr} = \frac{\omega^3}{9}\sqrt{\frac{2 \epsilon}{C (\widehat{C}_1^2 + \widehat{C}_2^2)}} \;.
\label{eq_2dof37}
\end{equation}

\begin{figure}[th]
\subfloat{\includegraphics[width=.45\textwidth]{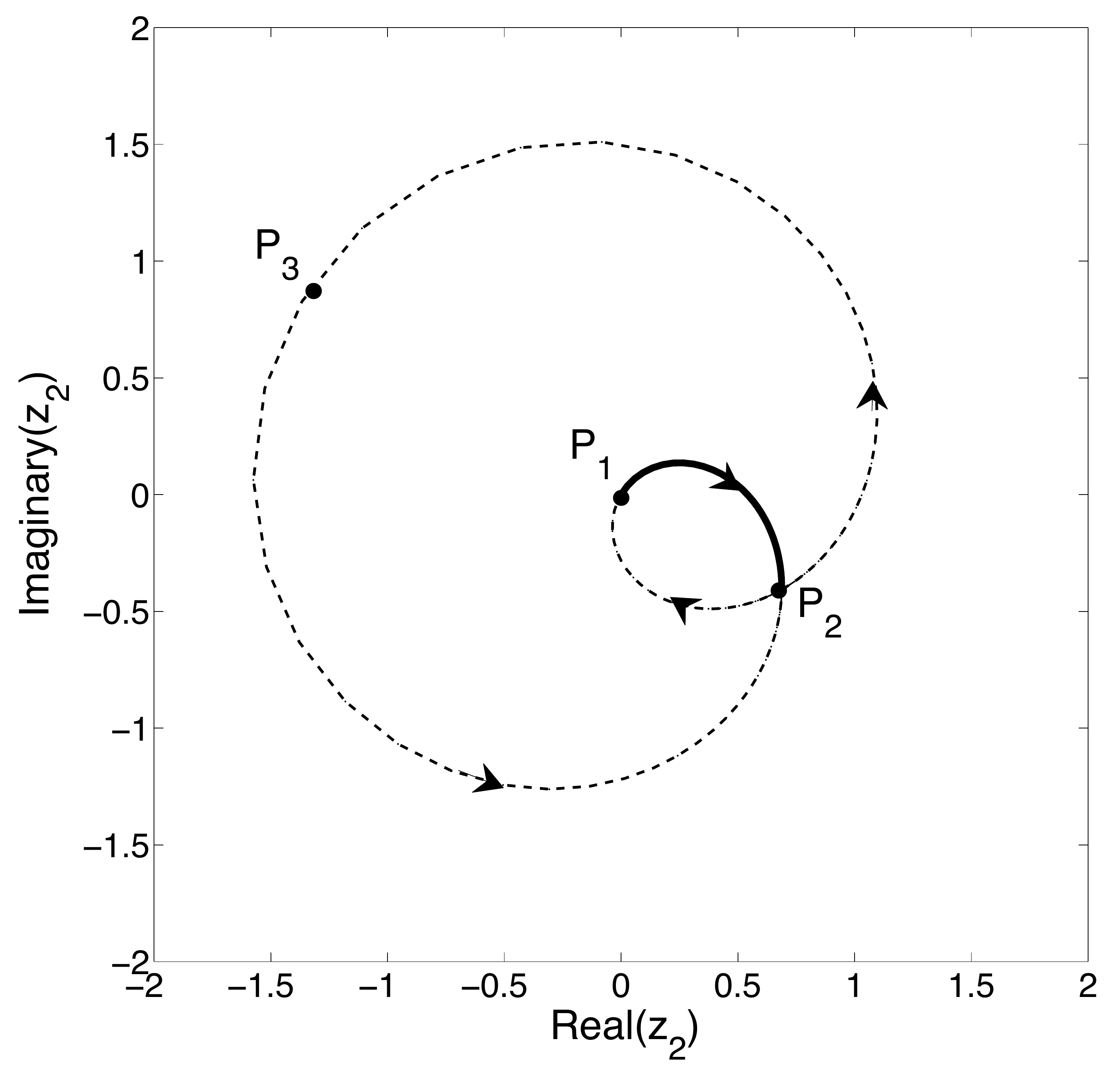}}
\hspace{.25in}
\subfloat{\includegraphics[width=.55\textwidth,height=.25\textwidth]{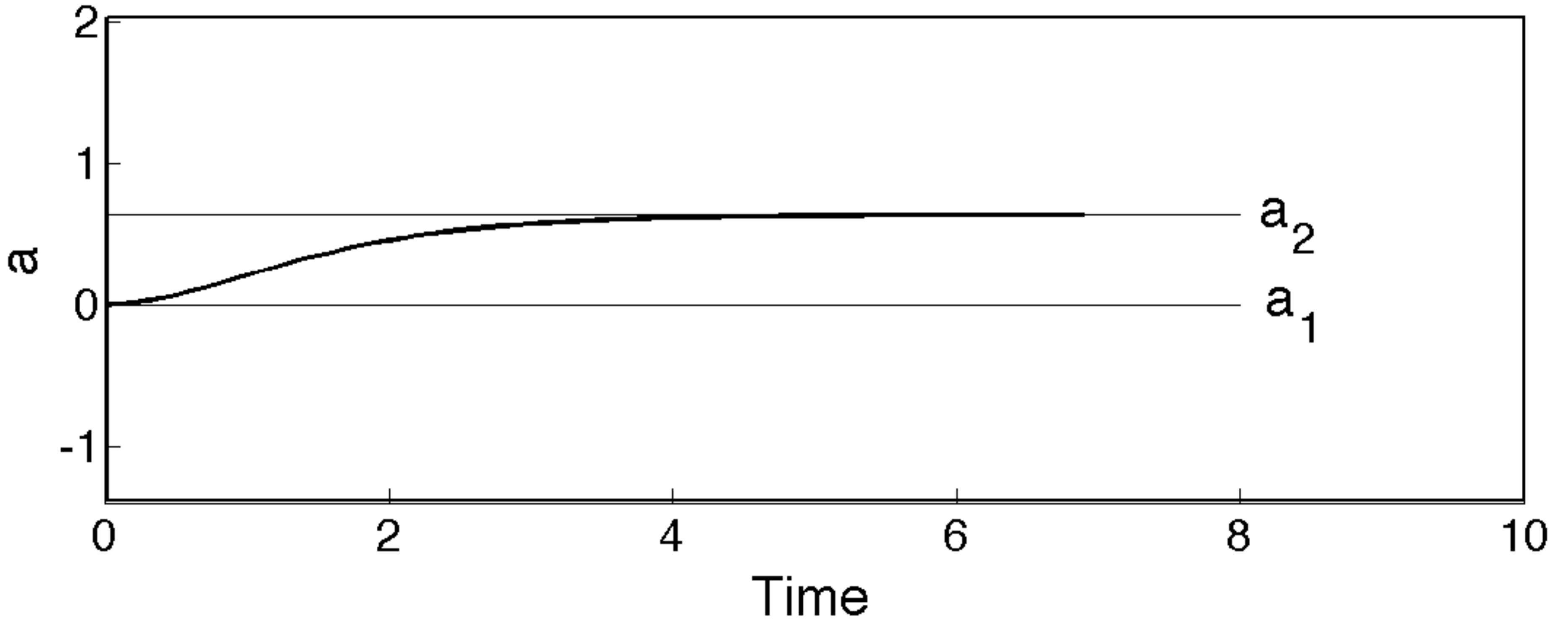}}
\caption{ {(a) Phase portrait of system given by Eq. (\ref{eq_2do31}) for an $z_3(0)=z_3(0)_{cr}$ given by Eq. (\ref{eq_2dof36}).\label{fig_2dof15}
	For this critical case, the system gets transferred from $P_1$ to $P_2$ via the trajectory shown in bold. The inner homoclinic loop also consists of the backward time dashed trajectory. The outer homoclinic loop is also shown. (b) The trajectory from $P_1$ to $P_2$ as function of time.} }%
\end{figure}

For this critical case, the three fixed points $P_i$ of the system given by Eq. (\ref{eq_2do31}) or Eq. (\ref{eq_2dof33}) have magnitudes $a_1=0,a_2=0.638, \text{ and } a_3=2.552$ respectively. In Fig. \ref{fig_2dof15}(a), the homoclinic orbit of the system with initial condition starting very close to $P_1$ is shown in bold. We confirm numerically that this orbit tends to $P_2$ as $t\rightarrow\infty$. The homoclinic loop can be completed by integrating backward in time from $P_1$ ($t\rightarrow-\infty$), and is shown as dashed orbit. There is another (`outer') homoclinic loop, also shown as dashed orbit. This loop also corresponds to the system given by Eq. (\ref{eq_2dof33}), but with different initial conditions (i.e. $z_2\neq0$), and hence is not important for our analysis. This loop starts at $P_2$ and approaches $P_3$ as $t\rightarrow\pm\infty$. The orbit from $P_1$ to $P_2$ is also shown as a function of time in Fig. \ref{fig_2dof15}(b).

We note that the expression for optimal velocity in Eq. (\ref{eq_2dof37}) is not a result of a formal optimization procedure. Such an optimization would most likely be analytically intractable, and computationally intensive because of the essential nonlinearity in the system. Our premise is that the analysis and computations carried out in this section provide insight into the energy transfer in the system. The resulting analytical criterion enables us to choose parameters that lead to capture into a homoclinic orbit, resulting in efficient energy transfer from the primary system to the attachment.  

\section{Tuned-Mass Damper: Energy Dissipation for Impulsive Excitation}\label{tmd}

\begin{figure}[ptbh]
\centering\includegraphics[width=0.3\textwidth]{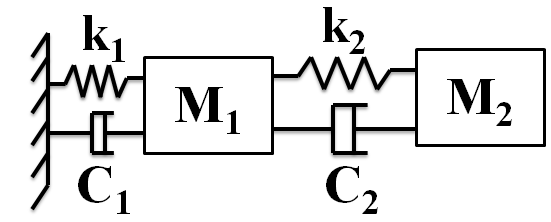}\caption{1-DOF
model for the tuned-mass damper}\label{fig:TMD}%
\end{figure}

Fig.  \ref{fig:TMD} shows a one-degree-of-freedom primary system with an added Tuned-Mass
Damper (TMD). The excitation is assumed to be impulsive, i.e. $\dot{x}%
_{1}\left(  0\right)  =v_{10}$. 

The equations of motion are given by%

\begin{gather}
m_{1}\ddot{x}_{1}+c_{1}\dot{x}_{1}+k_{1}x_{1}+c_{2}(\dot{x}_{1}-\dot{x}%
_{2})+k_{2}(x_{1}-x_{2})=0,\label{eq1}\\
m_{2}\ddot{x}_{2}+c_{2}(\dot{x}_{2}-\dot{x}_{1})+k_{2}(x_{2}-x_{1})=0.
\label{eq2}%
\end{gather}
With the introduction of
\begin{align}
\varepsilon &  =\frac{m_{2}}{m_{1}},\\
\kappa &  =\frac{k_{2}}{k_{1}},\\
\zeta_{1}  &  =\frac{c_{1}}{2\sqrt{m_{1}k_{1}}},\\
\zeta_{2}  &  =\frac{c_{2}}{2\sqrt{m_{1}k_{1}}},
\end{align}
we write Equations (\ref{eq1}, \ref{eq2}) in the following nondimensional form%
\begin{gather}
\ddot{x}_{1}+2\zeta_{1}\dot{x}_{1}+x_{1}+2\zeta_{2}(\dot{x}_{1}-\dot{x}%
_{2})+\kappa(x_{1}-x_{2})=0,\\
\varepsilon\ddot{x}_{2}+2\zeta_{2}(\dot{x}_{2}-\dot{x}_{1})+\kappa(x_{2}%
-x_{1})=0.
\end{gather}
The initial conditions are given by%
\begin{equation}
x_{1}\left(  0\right)  =0,\quad\dot{x}_{1}\left(  0\right)  =v_{10},\quad
x_{2}\left(  0\right)  =0,\quad\dot{x}_{2}\left(  0\right)  =0.
\end{equation}
In matrix form%
\begin{equation}
\mathrm{M}\mathbf{\ddot{x}}+\mathrm{C}\mathbf{\dot{x}}+\mathrm{K}\mathbf{x}=0,
\label{MCK}%
\end{equation}
where%
\begin{equation}
\mathbf{x=}\left(
\begin{array}
[c]{c}%
x_{1}\\
x_{2}%
\end{array}
\right)
\end{equation}
and%
\begin{equation}
\mathrm{M}=\left(
\begin{array}
[c]{cc}%
1 & 0\\
0 & \varepsilon
\end{array}
\right)  ,\quad\mathrm{C}=2\left(
\begin{array}
[c]{cc}%
\zeta_{1}+\zeta_{2} & -\zeta_{2}\\
-\zeta_{2} & \zeta_{2}%
\end{array}
\right)  ,\quad\mathrm{K}=\left(
\begin{array}
[c]{cc}%
1+\kappa & -\kappa\\
-\kappa & \kappa
\end{array}
\right)  .
\end{equation}
Now we recast Equation (\ref{MCK}) in first order form%
\begin{equation}
\mathbf{\dot{q}}=\mathrm{A}\mathbf{q}, \label{statespace}%
\end{equation}
where%
\begin{equation}
\mathbf{q}=\left(
\begin{array}
[c]{c}%
\mathbf{x}\\
\mathbf{\dot{x}}%
\end{array}
\right)  =\left(
\begin{array}
[c]{c}%
x_{1}\\
x_{2}\\
\dot{x}_{1}\\
\dot{x}_{2}%
\end{array}
\right)  ,
\end{equation}
and%
\begin{equation}
\mathrm{A}=\left(
\begin{array}
[c]{cc}%
0 & \mathrm{I}\\
-\mathrm{M}^{-1}\mathrm{K} & -\mathrm{M}^{-1}\mathrm{C}%
\end{array}
\right)  =\left(
\begin{array}
[c]{cccc}%
0 & 0 & 1 & 0\\
0 & 0 & 0 & 1\\
-1-\kappa & \kappa & -2\left(  \zeta_{1}+\zeta_{2}\right)  & 2\zeta_{2}\\
\frac{\kappa}{\varepsilon} & -\frac{\kappa}{\varepsilon} & 2\frac{\zeta_{2}%
}{\varepsilon} & -2\frac{\zeta_{2}}{\varepsilon}%
\end{array}
\right)  .
\end{equation}
The initial condition for the first order system (\ref{statespace}) is%
\begin{equation}
\mathbf{q}\left(  0\right)  =\left(
\begin{array}
[c]{c}%
0\\
0\\
v_{10}\\
0
\end{array}
\right)  . \label{IC}%
\end{equation}

\subsection{Optimal Parameters for Energy Dissipation}

In this Section, similarly to Wang et al. \cite{wang1984transient} we use
Lyapunov's second method to minimize an integral square performance measure of
damped vibrating structures subject to initial impulse. We aim to maximize the
quadratic cost function (energy damped by the tuned mass damper normalized with total
initial energy)%
\begin{equation}
\mathrm{J}=\frac{2\zeta_{2}%
{\displaystyle\int\limits_{0}^{\infty}}
(\dot{x}_{2}\left(  \tau\right)  -\dot{x}_{1}\left(  \tau\right)  )^{2}d\tau
}{\frac{1}{2}v_{10}^{2}}.
\end{equation}
This cost function can be expressed in terms of the state $\mathbf{q}$ as
\begin{gather}
\mathrm{J}=\frac{4\zeta_{2}}{v_{10}^{2}}%
{\displaystyle\int\limits_{0}^{\infty}}
\mathbf{q}^{\mathrm{T}}\left(  \tau\right)  \mathrm{Q}\mathbf{q}\left(
\tau\right)  d\tau,\\
\mathrm{Q}=\frac{4\zeta_{2}}{v_{10}^{2}}\left(
\begin{array}
[c]{cccc}%
0 & 0 & 0 & 0\\
0 & 0 & 0 & 0\\
0 & 0 & 1 & -1\\
0 & 0 & -1 & 1
\end{array}
\right)  .
\end{gather}

The matrix $\mathrm{Q}$ is symmetric, positive semidefinite. According to
Lyapunov theory (Lyapunov's second method \cite{khalil1996nonlinear}), for a
stable system Eq. (\ref{statespace}) there exists a positive semidefinite
$\mathrm{P}$ satisfying the (Lyapunov) equation%
\begin{equation}
-\mathrm{Q}=\mathrm{A}^{\mathrm{T}}\mathrm{P}+\mathrm{PA.} \label{Lyapunov}%
\end{equation}
Left multiplying with $\mathbf{q}^{\mathrm{T}}$ and right multiplying with
$\mathbf{q}$ yields
\begin{equation}
-\mathbf{q}^{\mathrm{T}}\mathrm{Q}\mathbf{q}=\mathbf{q}^{\mathrm{T}%
}\mathrm{A^{\mathrm{T}}P}\mathbf{q}+\mathbf{q}^{\mathrm{T}}\mathrm{PA}%
\mathbf{q.}%
\end{equation}
Using Eq. (\ref{statespace}) we get%
\begin{equation}
-\mathbf{q}^{\mathrm{T}}\mathrm{Q}\mathbf{q}=\mathbf{\dot{q}}^{\mathrm{T}%
}\mathrm{P}\mathbf{q}+\mathbf{q}^{\mathrm{T}}\mathrm{P}\mathbf{\dot{q}}%
=\frac{d}{dt}\left(  \mathbf{q}^{\mathrm{T}}\mathrm{P}\mathbf{q}\right)  ,
\end{equation}
and thus%
\begin{equation}
\mathrm{J}=%
{\displaystyle\int\limits_{0}^{\infty}}
-\frac{d}{dt}\left(  \mathbf{q}^{\mathrm{T}}\mathrm{P}\mathbf{q}\right)
d\tau=\mathbf{q}^{\mathrm{T}}\left(  0\right)  \mathrm{P}\mathbf{q}\left(
0\right)  -\mathbf{q}^{\mathrm{T}}\left(  \infty\right)  \mathrm{P}%
\mathbf{q}\left(  \infty\right)  .
\end{equation}
The system Eq. (\ref{statespace}) is asymptotically stable, so $\mathbf{q}\left(
\infty\right)  \rightarrow\mathbf{0}$, and%
\begin{equation}
\mathrm{J}=\mathbf{q}^{\mathrm{T}}\left(  0\right)  \mathrm{P}\mathbf{q}%
\left(  0\right)  =v_{10}^{2}\mathrm{P}_{33}. \label{J}%
\end{equation}
Solving the Lyapunov equation Eq. (\ref{Lyapunov}) for $\mathrm{P}$ and
substituting into Eq. (\ref{J}) yields%
\begin{equation}
\mathrm{J}=\frac{\text{$\zeta_{2}$}\varepsilon\left(  4\text{$\zeta_{1}%
^{2}\zeta_{2}$}\kappa+\text{$\zeta_{2}$}(4\text{$\zeta$}_{1}\text{$\zeta_{2}$%
}+\varepsilon)+\text{$\zeta$}_{1}\kappa^{2}(\varepsilon+1)\right)
}{\text{$\zeta$}_{1}\text{$\zeta_{2}$}\left(  4\text{$\zeta$}_{1}%
\text{$\zeta_{2}$}\kappa+4\text{$\zeta_{2}$}^{2}+\kappa^{2}\right)
+\varepsilon^{2}(\text{$\zeta$}_{1}+\text{$\zeta_{2}$})\left(  \text{$\zeta$%
}_{1}\kappa^{2}+\text{$\zeta_{2}$}\right)  +2\text{$\zeta$}_{1}\text{$\zeta
_{2}$}\varepsilon\left(  \kappa(2\text{$\zeta$}_{1}(\text{$\zeta$}%
_{1}+\text{$\zeta_{2}$})-1)+2\text{$\zeta_{2}$}(\text{$\zeta$}_{1}%
+\text{$\zeta_{2}$})+\kappa^{2}\right)  }.
\end{equation}
Solving $\frac{\partial J}{\partial\kappa}=0$ yields%
\begin{equation}
\kappa=\frac{\varepsilon+2\text{$\zeta$}_{1}\text{$\zeta_{2}$}}{1+\varepsilon
-2\text{$\zeta$}_{1}^{2}}. \label{kappa}%
\end{equation}
Substituting Eq. (\ref{kappa}) into $\frac{\partial J}{\partial\text{$\zeta_{2}$}%
}=0$ yields%
\begin{equation}
(\varepsilon+1)\varepsilon^{3}+8\text{$\zeta$}_{1}\text{$\zeta_{2}$%
}\varepsilon^{2}\left(  \varepsilon+1-\text{$\zeta$}_{1}^{2}\right)
+4\text{$\zeta_{2}^{2}$}\left(  \varepsilon+1-\text{$\zeta$}_{1}^{2}\right)
\left(  4\text{$\zeta$}_{1}^{2}(\varepsilon+1-\text{$\zeta$}_{1}%
^{2})-\varepsilon-1\right)  =0.
\end{equation}
Solving for $\zeta_{2}$ results in%
\begin{equation}
\text{$\zeta_{2}$}=\frac{2\varepsilon^{2}\zeta_{1}\left(  1+\varepsilon
-\zeta_{1}^{2}\right)  +\varepsilon\left\vert 1+\varepsilon-2\zeta_{1}%
^{2}\right\vert \sqrt{\varepsilon\left(  1+\varepsilon-\zeta_{1}^{2}\right)
}}{2\left(  1+\varepsilon-\zeta_{1}^{2}\right)  \left(  1+\varepsilon-4\left(
1+\varepsilon\right)  \zeta_{1}^{2}+4\zeta_{1}^{4}\right)  }.\label{Eq:zeta2_opt}
\end{equation}

\begin{figure}[ptbh]
\centering\includegraphics[width=0.8\textwidth]{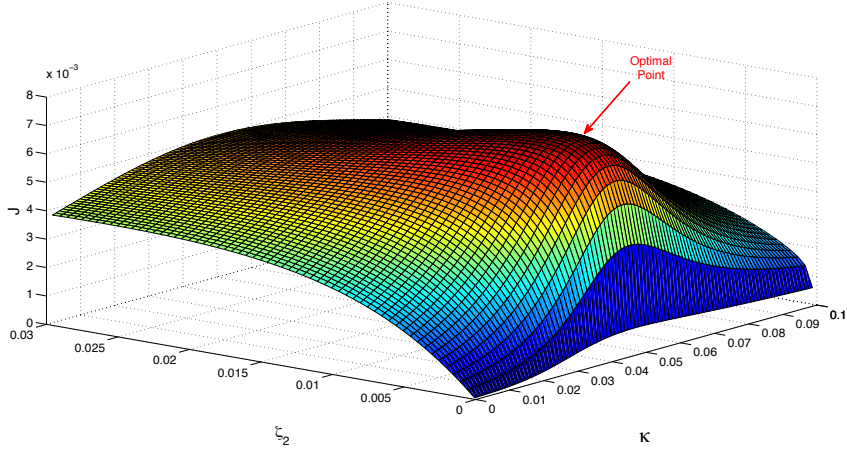}
\caption{1 DOF primary system with TMD:Energy dissipation J$\left(  \zeta_{2},\kappa\right)  $ and point of
optimal dissipation ($\zeta_{2}=0.0055,\kappa=0.048$). The parameter values
are $\varepsilon=0.05$ and $\zeta_{1}=0.02$.}%
\label{fig:Jwithoptimalpoint}%
\end{figure}

Fig.  \ref{fig:Jwithoptimalpoint} shows the energy dissipation
$\mathrm{J}$ as a function of $\zeta_{2}$ and $\kappa$ for $\varepsilon=0.05$
and $\zeta_{1}=0.02$. The maximum (optimal) value is $\mathrm{J}%
=0.00726$. 

\begin{figure}[ptbh]
\centering\includegraphics[width=0.8\textwidth]{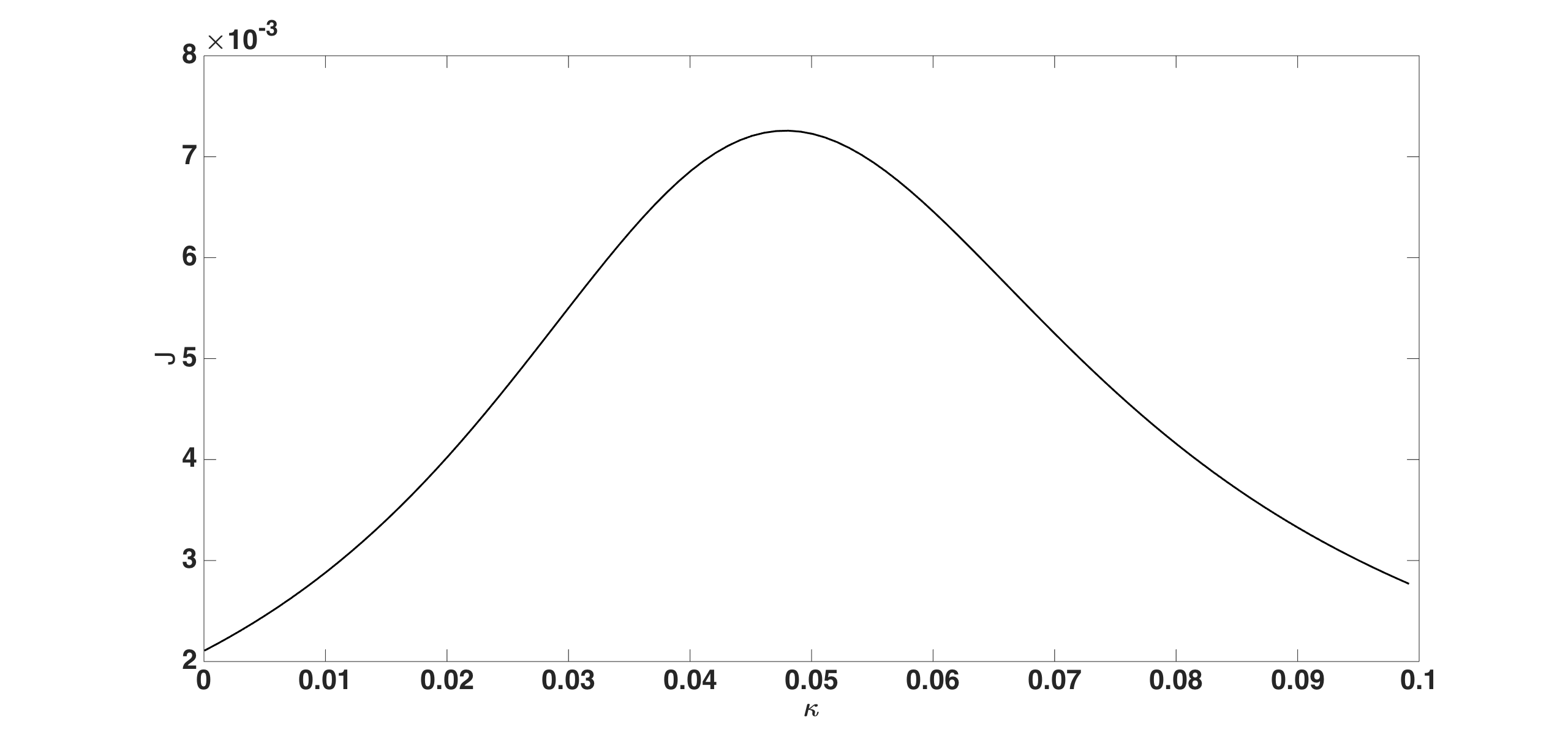} \caption{1 DOF primary system with TMD: Energy
dissipation \textrm{J}$\left(  \kappa\right)  $ at optimal $\zeta_{2}%
=0.0055$.}%
\label{fig:Jvskappa}%
\end{figure}

Fig.  \ref{fig:Jvskappa} shows the energy dissipation
$\mathrm{J}$ as a function of $\kappa$ for $\varepsilon=0.05$ and $\zeta
_{1}=0.02$ at the optimal $\zeta_{2}=0.0055$. This numerically computed optimal result matches exactly with the analytical value of $\zeta_2$ computed using Eq. (\ref{Eq:zeta2_opt}).

\subsection{Response of 2 DOF Primary System with TMD}
 The non-dimensional equations of motion for a two-degree-of-freedom primary system attached a with
a tuned mass damper (TMD), similar to the system shown in Fig.
\ref{fig_2dof01} can be written as
\begin{subequations}
\begin{gather}
x_{1}^{\prime\prime}+2\zeta_{1}x_{1}^{\prime}+2\zeta_{12}(x_{1}^{\prime}%
-x_{2}^{\prime})+x_{1}+k_{12}(x_{1}-x_{2})=0,\\
\mu x_{2}^{\prime\prime}+2\zeta_{12}(x_{2}^{\prime}-x_{1}^{\prime})+2\zeta
_{3}(x_{2}^{\prime}-x_{3}^{\prime})+k_{12}(x_{2}-x_{1})+k_{tmd}(x_{2}%
-x_{3})=0,\\
\epsilon x_{3}^{\prime\prime}+2\zeta_{3}(x_{3}^{\prime}-x_{2}^{\prime
})+k_{tmd}(x_{3}-x_{2})=0.
\end{gather}\label{eq_2doftmd}
\end{subequations}

Similar to the 1 degree-of-freedom case, Lyapunov's second method can be used to maximize an integral square performance measure subject to initial impulse.
Given all other parameters in the above equations, we want to choose parameters $(\zeta_3,k_{tmd})$ to minimize J, the total normalized energy dissipation in the TMD
\begin{equation}
J=\dfrac{2\zeta_3\int_{0}^{\infty} (\dot{x_2}(\tau)-\dot{x_3}(\tau))^2d\tau}{\frac{1}{2}v^2_{10}}
\end{equation}

\begin{figure}[ptbh]
\centering\includegraphics[width=0.8\textwidth]{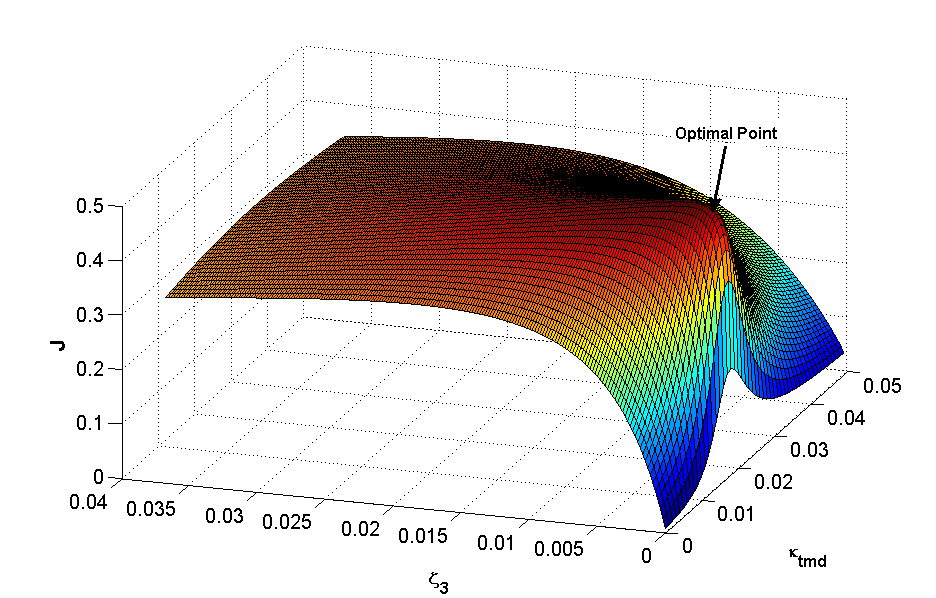}
\caption{2 DOF primary with TMD: Energy dissipation J$\left(  \zeta_{3},k_{tmd}\right)  $ and point of
optimal dissipation ($\zeta_{3}=0.0039,k_{tmd}=0.02$). The parameter values
are $\varepsilon=0.0318$, $\mu=0.6364$, $\zeta_1=.0074$, $k_{12}=2.5$ and $\zeta_{12}=0.0148$.}%
\label{fig:2dofJwithoptimalpoint}%
\end{figure}

\begin{figure}[th]
\centering
\includegraphics[width=0.6\textwidth]{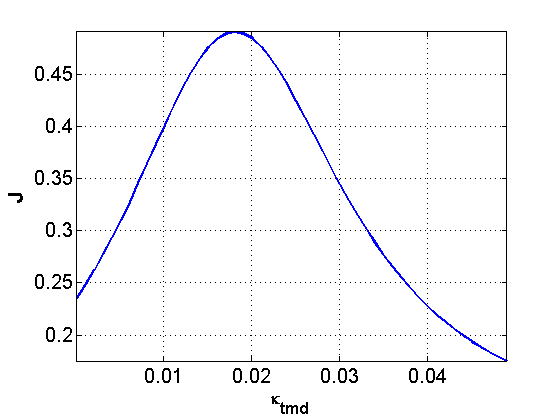}
\caption{2 DOF primary with TMD: Quadratic cost function $J$ for different values of TMD stiffness $k_{tmd}$ at optimal damping $\zeta_3=.0039$.}%
\label{J_2dof_tmd}
\end{figure}

\begin{figure}[th]
\centering
\includegraphics[scale=0.5]{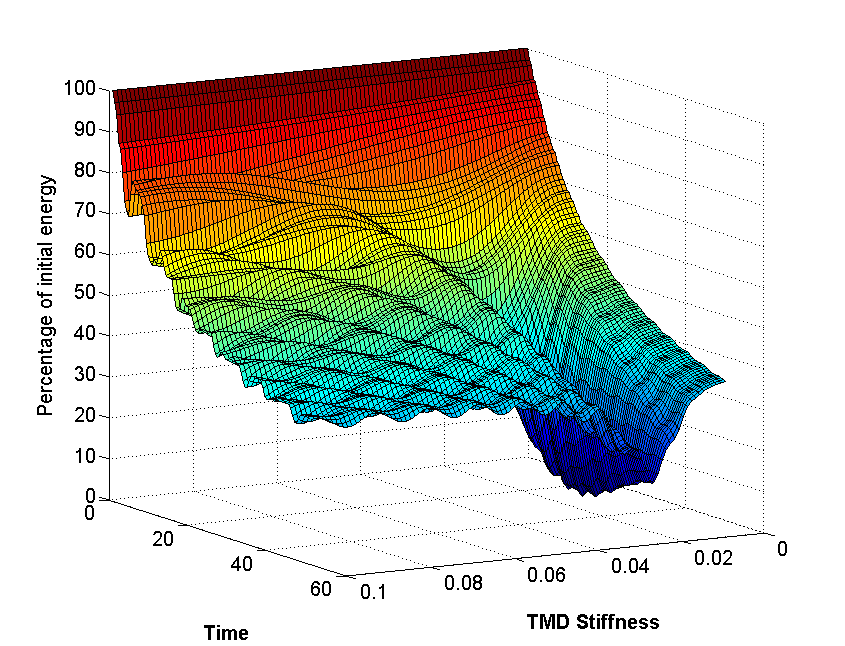}\caption{Response of the 2 DOF system
with the TMD for different values of TMD stiffness.}\label{fig_2dofnestmd02}%
\end{figure}

 An analytical expression for the optimal value quadratic cost function $J$, and corresponding parameter values in the two-degree-of-freedom case could not be determined. The Lyapunov equation in this case was solved numerically, and a gradient ascent algorithm was applied to find the optimal parameter values. Fig.  \ref{fig:2dofJwithoptimalpoint} shows the energy dissipation in TMD as function of $\zeta_3$ and $k_{tmd}$, along with the optimal point. Fig.  \ref{J_2dof_tmd} shows the variation in the cost function with change in TMD stiffness. The response of the system given in Eq. (\ref{eq_2doftmd}) for range of $k_{tmd}$ values for a particular initial impulse velocity given to mass $M_{1}$ is given in Fig. (\ref{fig_2dofnestmd02}). As can be observed from this figure, best TMD performance is achieved at $k_{tmd} = 0.0181$ . This value is close to the optimal value found by Lyapunov analysis, and hence confirms the validity of our analysis. The optimal value is also independent of magnitude of the initial velocity given to mass $M_{1}$, since the problem is linear.
 
\clearpage

\section{Results}\label{results}
\subsection{NES Response}

\begin{figure}[th]
\centering
\includegraphics[scale=0.5]{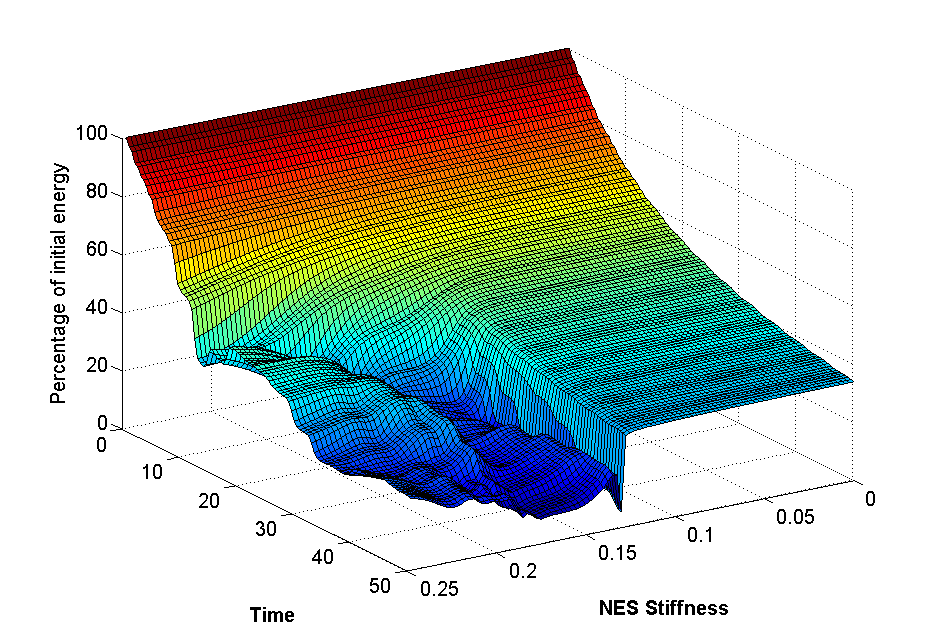}\caption{Response of the 2 DOF system
with the NES for different values of NES stiffness and an initial velocity of
$v_{0}=0.2$ .}\label{fig_2dofnestmd03}%
\end{figure}

Unlike the TMD, the NES response depends on the magnitude of the initial
impulse velocity given to mass $M_{1}$. Therefore, for comparison sake, if it
is assumed that the aim of the vibration mitigation mechanism is to attenuate
the vibrations caused by a particular velocity then, the performance of the
NES for different stiffness ($C$) can be compared with that of the TMD. Fig. 
\ref{fig_2dofnestmd03} shows the response of the NES to a particular initial
velocity ($v_{0}=0.2$) for a range of NES stiffness. 

\subsection{Response of TMD and NES to Variations in Parameters}

\begin{figure}[th]
\centering
\par
\includegraphics[scale=0.5]{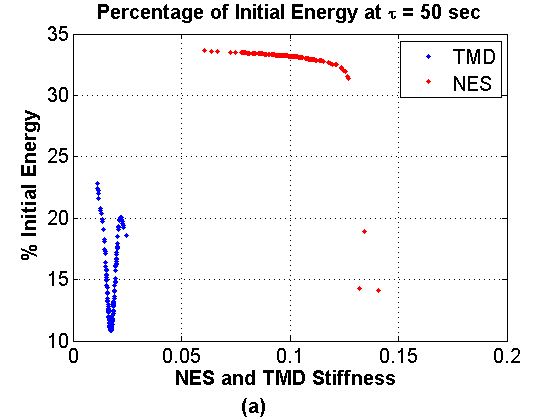}
\includegraphics[scale=0.5]{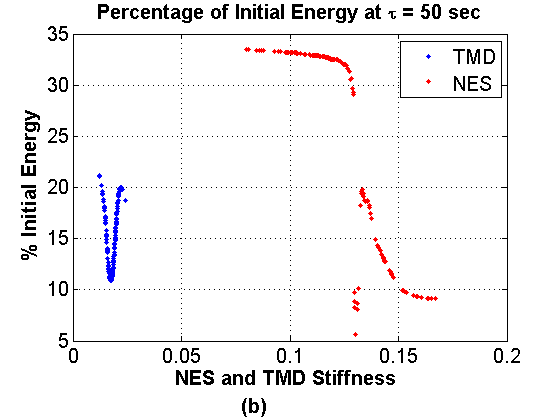} \caption{Comparison of the
performance of TMD and the NES with a normal distribution of stiffness. Part
(a) uses Eq. (\ref{eq_2dofvcr}) exactly to obtain the NES stiffness, whereas
part in (b), the design critical velocity is 90 $\%$ of the actual initial
impulse velocity.}%
\label{fig_2dofnestmd04}%
\end{figure}

Robustness is an important consideration in the design of either the TMD or
NES system. To get more insight into the performance of the TMD or NES to
changes in system parameters, several simulations have been carried out. For
the design of NES systems, the following relationship (derived in section 4) between the critical
initial velocity and system parameters is used,
\begin{equation}
v_{cr}=\frac{\omega^{3}}{9}\sqrt{\frac{2\epsilon}{C(\widehat{C}_{1}%
^{2}+\widehat{C}_{2}^{2})}}, \label{eq_2dofvcr}%
\end{equation}
Using the Eq. (\ref{eq_2dofvcr}), the value of the NES stiffness $C$ can be
calculated for a particular initial velocity. For the system given in Eq. (\ref{eq_2dof01}), this value of critical velocity was 0.09 which is close to the numerically determined value of 0.115. From Fig. \ref{fig_2dofnestmd02} and
\ref{fig_2dofnestmd03}, it is clear that both TMD and NES have an optimal
stiffness (for the NES case for a particular initial velocity of mass $M_{1}%
$). This allows for the selection of a nominal optimal value. Fig. 
\ref{fig_2dofnestmd04} (a) and (b) show the results of a simulation for
different NES stiffnesses. In Fig. \ref{fig_2dofnestmd04} (a), the NES
stiffness is designed using Eq. (\ref{eq_2dofvcr}) and then the simulations
are done by assuming that the actual stiffness has a normal probability
distribution function (pdf) about the nominal optimal value thus calculated.
In Fig. \ref{fig_2dofnestmd04} (b), the NES stiffness is designed for 90 $\%$
of the desired critical velocity. This is done because as can be concluded
from part (a) the relationship between the critical velocity and the NES
stiffness is not exactly accurate as it is based on several layers of
approximations. The stiffness value for the NES is again assumed to have a
normal pdf around the new nominal value thus calculated (using critical
velocity which is 90 $\%$ of the critical velocity) and the results of the
simulation are presented in Fig. \ref{fig_2dofnestmd04} (b). In both Fig.
\ref{fig_2dofnestmd04} (a) and (b), the stiffness for the TMD is assumed to
have a normal pdf with the mean equal to the optimum value obtained using computations of Section \ref{tmd}.

\begin{figure}[th]
\centering
\includegraphics[scale=0.5]{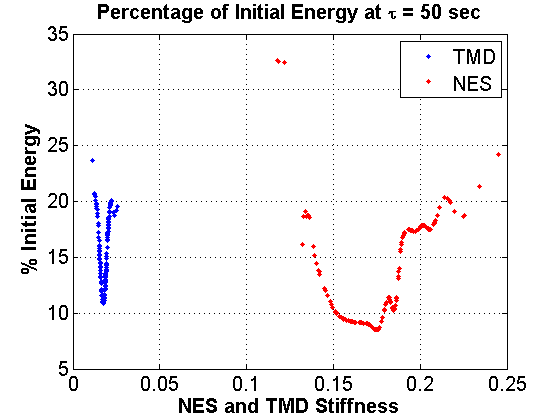} \caption{Comparison of the
performance of TMD and the NES with a normal distribution of Stiffnesses. The
design critical velocity is 75 $\%$ of the actual initial impulse velocity.}%
\label{fig_2dofnestmd05}%
\end{figure}

Fig.  \ref{fig_2dofnestmd05} contains the results for the designed critical
velocity being 75$\%$ of the actual initial impulse velocity.

One important aspect of NES performance is its dependence on initial velocity.
Also, it has been previously observed (\cite{sapsis_01}) that the performance
of NES increases (or degrades) significantly if the initial impulse velocity
given to mass $M_{1}$ is above (or below) a certain threshold. From
performance perspective, it is desirable to be as close to this critical limit
as possible, however, such a setup will have low robustness as even a little
change in the NES stiffness can cause significant degradation in performance.
An alternative design strategy could be to keep the designed critical velocity
lower than the initial impulse velocity given to mass $M_{1}$. Next, we study the effect
of system parameter variation on TMD or NES performance.

\subsection{Variation in Mass $M_{2}$}

\begin{figure}[th]
\centering
\includegraphics[width=.45\textwidth]{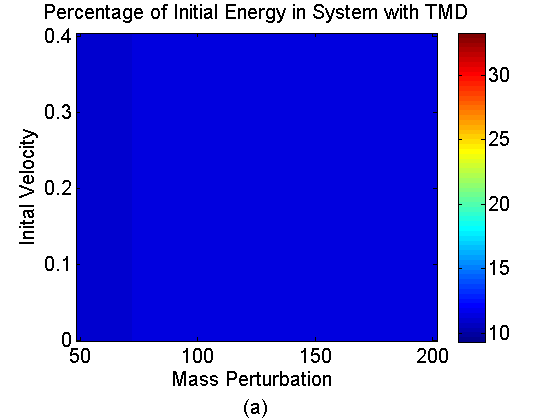}
\includegraphics[width=.4\textwidth]{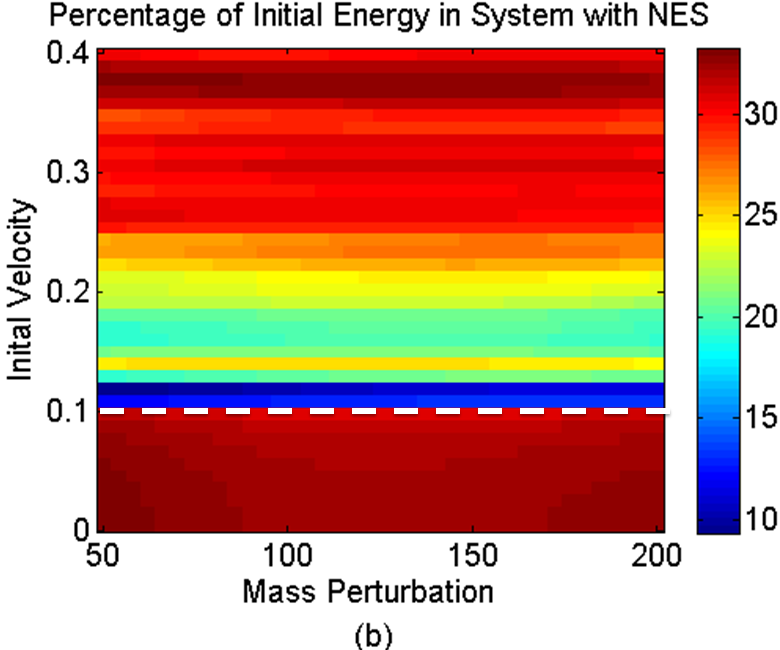} \caption{Comparison of the
performance of TMD and the NES with perturbations in the value of mass $M_{2}$
(Results at $\tau= 50$). The NES was designed to keep $v_{0} = 0.1$ as the
critical velocity (shown as white dashed line).}%
\label{fig_2dofnestmd06}%
\end{figure}

The performance of the TMD is compared with the performance of the NES for
perturbations in the mass $M_{2}$ in Fig. \ref{fig_2dofnestmd06}. All
other system parameters are kept constant. The sensitivity of TMD performance to mass perturbation is
found to be minimal. For the NES, the performance is highly dependent on the initial
velocity. In this particular simulation the NES stiffness was designed to keep
$v_{0} = 0.1$ as the critical velocity (by putting $v_{cr} = 0.09$ in Eq.
(\ref{eq_2dofvcr})). The NES however is able to maintain its performance over
the whole range of mass $M_{2}$ perturbations shown in  Fig. \ref{fig_2dofnestmd06}.

\subsection{Variation in Mass $M_{3}$}

\begin{figure}[th]
\centering
\includegraphics[width=.45\textwidth]{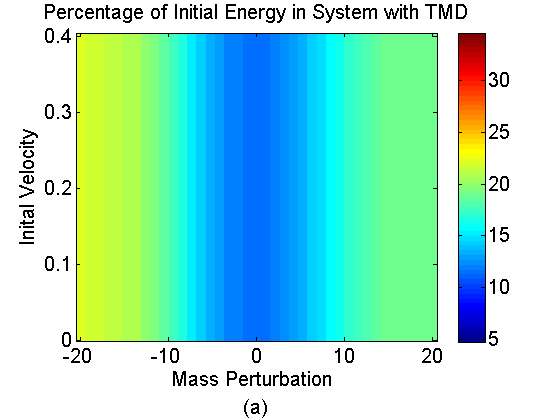}
\includegraphics[width=.4\textwidth]{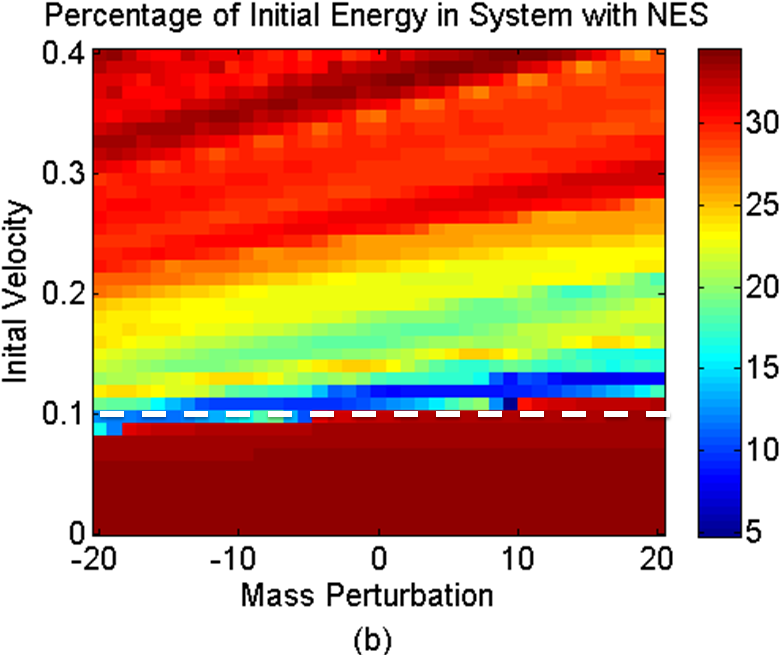} \caption{Comparison of the
performance of TMD and the NES with perturbations in the value of mass $M_{3}$
(Results at $\tau= 50$). The NES was designed to keep $v_{0} = 0.1$ as the
critical velocity (shown as white dashed line).}%
\label{fig_2dofnestmd07}%
\end{figure}

The performance of the TMD compared with the performance of the NES for
perturbations in the mass $M_{3}$ is shown in Fig. \ref{fig_2dofnestmd07}. All
other system parameters are kept constant and only the mass $M_{3}$ is varied around its nominal value. As can be observed from Fig.
\ref{fig_2dofnestmd07}, the TMD performance is more sensitive to
changes in mass $M_{3}$ as compared to the mass $M_{2}$ , though this may be
explained partially by the much lower nominal value of $M_{3}$ as compared to
$M_{2}$. The NES shows robustness to the changes in mass $M_{3}$. However the dependence of NES performance on initial velocity makes
its performance bad in some velocity regimes.

\subsection{Variation in Damping Coefficient $\zeta_{3}$}

\begin{figure}[th]
\centering
\includegraphics[width=.45\textwidth]{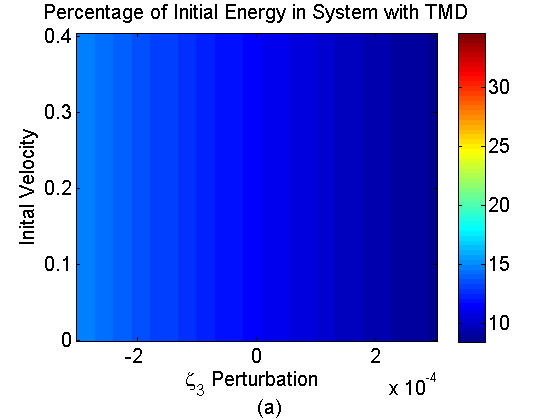}
\includegraphics[width=.4\textwidth]{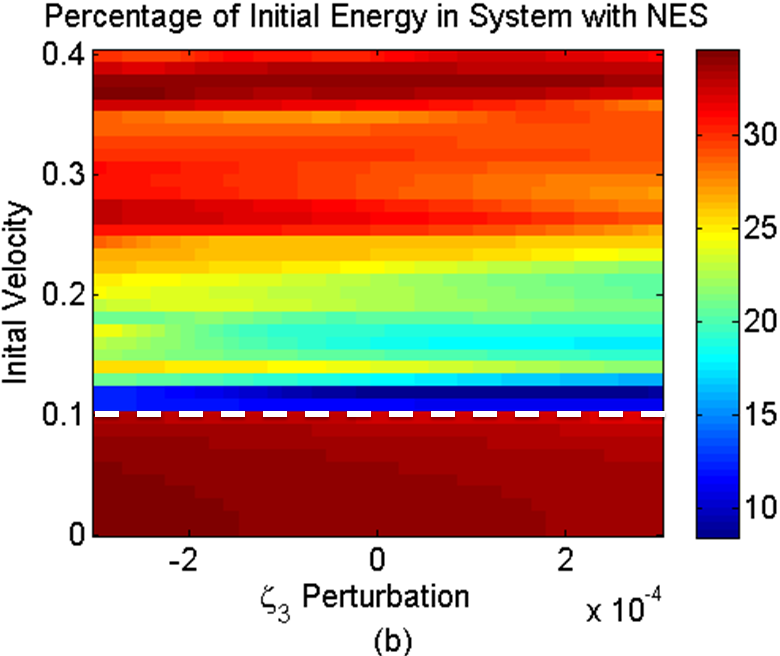} \caption{Comparison of the
performance of TMD and the NES with perturbations in the value of damping
coefficient $\zeta_{3}$ (Results at $\tau= 50$). The NES was designed to keep
$v_{0} = 0.1$ as the critical velocity (shown as white dashed line).}%
\label{fig_2dofnestmd08}%
\end{figure}

The performance of the TMD compared with the performance of the NES for
perturbations in the damping coefficient $\zeta_{3}$ is shown in Fig.
\ref{fig_2dofnestmd08}. All other system parameters are kept constant and only
the damping coefficient $\zeta_{3}$ is varied for a small range around its
nominal value. Damping coefficients in dampers may change with use and thus
this simulation can provide an idea of the change in TMD or NES performance to
change in the damping coefficient value $\zeta_{3}$. As can be observed from
Fig. \ref{fig_2dofnestmd08}, the performance of the NES to the changes in
$\zeta_{3}$ is again more robust than TMD near the designed velocity. 
However, again the dependence of NES performance on initial velocity makes its performance bad in some velocity regimes.

\subsection{Comparison of TMD and NES Performance in Low Damping Scenario}

\begin{figure}[th]
\centering
\includegraphics[width=0.45\textwidth]{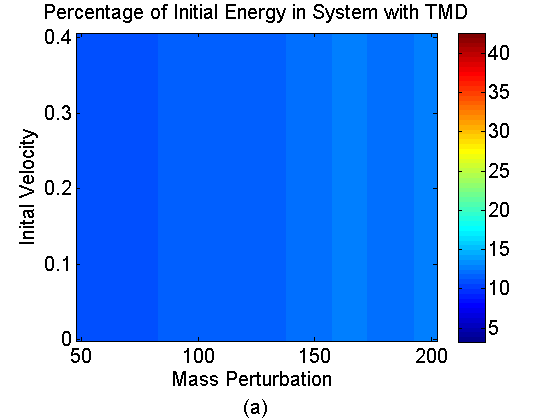}
\includegraphics[width=0.4\textwidth]{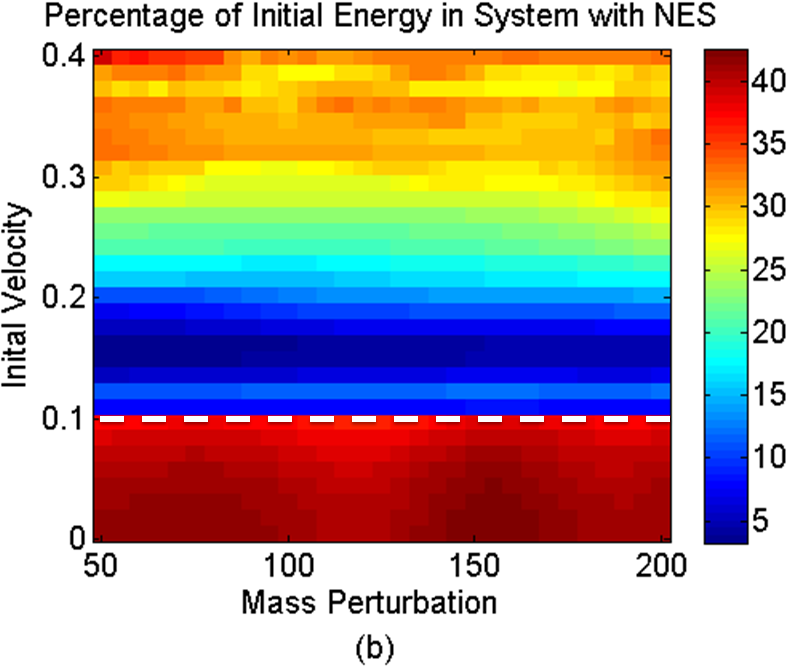} \caption{Comparison of the
performance of TMD and the NES with perturbations in the mass $M_{2}$. The
results are at $\tau=130$.The NES was designed to keep $v_{0} = 0.1$ as the critical velocity (shown as white dashed line).}%
\label{fig_ld_2dof_01}%
\end{figure}

\begin{figure}[th]
\centering
\includegraphics[width=0.45\textwidth]{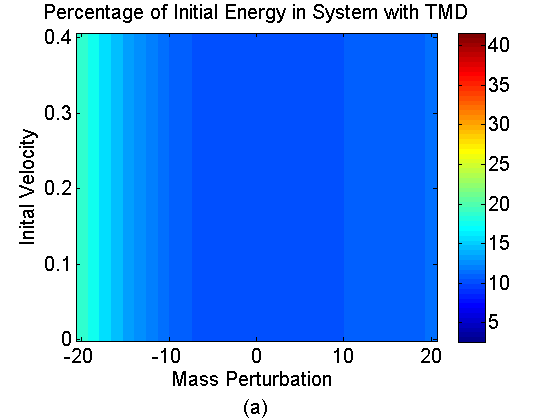}
\includegraphics[width=0.4\textwidth]{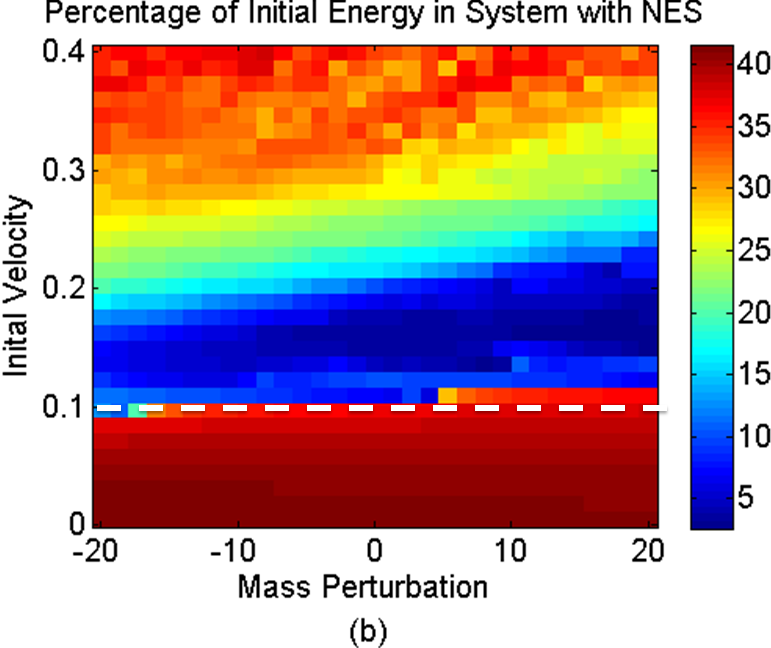} \caption{Comparison of the
performance of TMD and the NES with perturbations in the mass $M_{3}$. The
results are at $\tau= 130$. The NES was designed to keep
$v_{0} = 0.1$ as the critical velocity (shown as white dashed line).}%
\label{fig_ld_2dof_02}%
\end{figure}

\begin{figure}[th]
\centering
\includegraphics[width=0.45\textwidth]{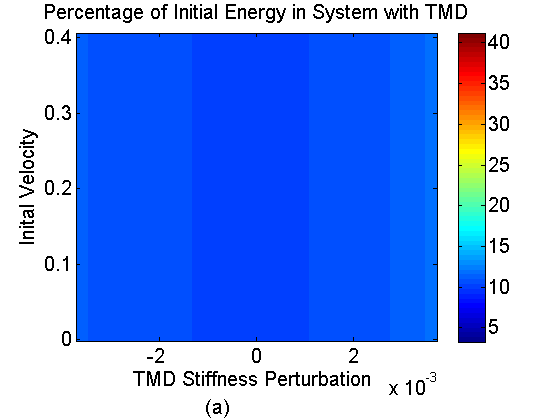}
\includegraphics[width=0.4\textwidth]{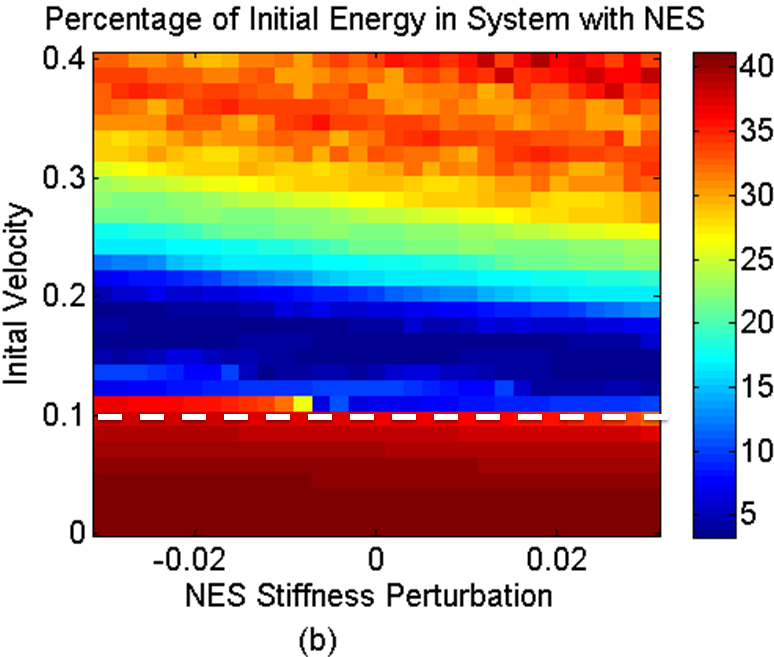} \caption{Comparison of the
performance of TMD and the NES with perturbations in the stiffness $k_{3}$ and
$C$. The results are at $\tau= 130$. The NES was designed to keep
$v_{0} = 0.1$ as the critical velocity (shown as white dashed line).} %
\label{fig_ld_2dof_03}%
\end{figure}

\begin{figure}[th]
\centering
\includegraphics[width=0.45\textwidth]{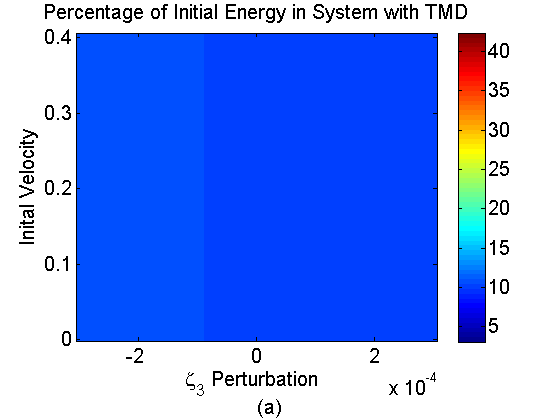}
\includegraphics[width=0.4\textwidth]{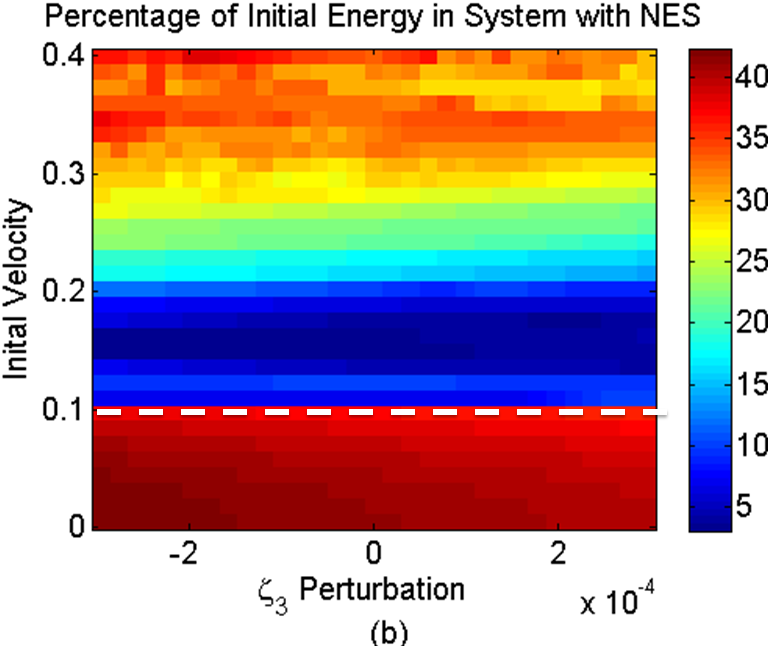} \caption{Comparison of the
performance of TMD and the NES with perturbations in damping ratio $\zeta_{3}%
$. The results are at $\tau= 130$. The NES was designed to keep
$v_{0} = 0.1$ as the critical velocity (shown as white dashed line). }%
\label{fig_ld_2dof_04}%
\end{figure}

Now we compare the performance of the NES and TMD for the case when primary damping is low compared to the nominal design considered earlier. The damping values considered for the results presented
in the subsequent part are $b_{1}=50$ Ns/m, $b_{2}=50$ Ns/m and $b_{3}=130$ Ns/m. The NES
performance is expected to improve for primary systems with lower damping (as reported in
\cite{sapsis_01}). The results for the simulations for variations in masses
$M_{2}$ and $M_{3}$ as well as the stiffnesses of TMD and NES and damping
ratio $\zeta_{3}$ are given in Figs. \ref{fig_ld_2dof_01}, \ref{fig_ld_2dof_02}%
, \ref{fig_ld_2dof_03} and \ref{fig_ld_2dof_04}. The time horizon is now increased to 130. It is clear from these figures that for the given time horizon and close to optimal input velocity, NES can outperform TMD for the case of no parametric perturbations (i.e. the best performance of NES is better than best performance of TMD). Close to the optimal input velocity, NES also shows better or similar robustness to parametric perturbations compared to TMD  in all four cases. Furthermore, the robustness to perturbations in input velocity is also increased compared to the case when primary damping was higher. 

\section{Conclusions}

Using a combination of analytical and numerical techniques, we have presented a framework for computing system parameters that lead to efficient one-way energy transfers in multi-degree-of-freedom systems with NES. We carry out explicit calculations for NES attached to a two-degree-of-freedom primary system, and the process of generalizing to N($\geq 2$) degree-of-freedom systems is also described. We exploit the separation of time-scales due to weak damping in the primary system, and implement complexification-averaging in the full system, followed by dimensionality reduction. Numerical evidence is presented for validity of this dimension reduction in the slow time scale. The analysis confirms that the near-optimal initial conditions for energy transfer to NES are close to homoclinic orbit in the undamped system on a slower time scale, as discovered earlier for the case of single-degree-of-freedom primary systems. Hence, our work generalizes previous work in parameter optimization of NES attached to single-degree-of-freedom systems, and provides evidence that homoclinic orbits of the underlying Hamiltonian system play a crucial role in efficient nonlinear energy transfers, even in high dimensional systems.

The performance of the optimally tuned NES and TMD attachments are compared under perturbations to primary system parameters. It is observed that NES performance is robust even though its best performance may not be as good as a perfectly tuned TMD. For the case of weakly damped primary system and a range of input velocities, the NES shows better robustness than the perfectly tuned TMD to parametric perturbation. Our comparisons also show degradation of NES performance far away from designed impulsive velocity, as has been previously observed \cite{vakphysD01}.

While our current work has considered NES with cubic essential nonlinearity of the spring, recent works \cite{al2013numerical, liu2014tailoring,gendelman2015dynamics,gourc2015quenching} have explored other nonlinear designs of both the spring and damper. Initial studies and numerical parameter optimization efforts have shown that these alternative nonlinear designs reduce the initial energy dependence of the cubic NES, and allow for a more graceful degradation in performance as initial energy moves away from the design point. Extensions of our work to such nonlinear attachments will be a topic of future work.
\section{Acknowledgements}
The authors would like to thank the reviewers for their careful reading of the manuscript and many constructive comments.
\appendix
\section{Analysis of an N-DOF System with NES}
\label{n-dof}

The analysis conducted on the two-degree-of-freedom NES can also be used to
study the vibration mitigation of a generic system consisting of $n$ ($> 2$)
masses connected in a series with a nonlinear energy sink (NES) attached at
the end. The equations of motion for such a system can be written as,%

\begin{subequations}
\begin{gather}
M_{1}\ddot{x}_{1}+b_{0}\dot{x}_{1}+b_{1}(\dot{x}_{1}-\dot{x}_{2})+\kappa
_{0}x_{1}+\kappa_{1}(x_{1}-x_{2})=0,\\
M_{i}\ddot{x}_{i}+b_{i-1}(\dot{x}_{i}-\dot{x}_{i-1})+b_{i}(\dot{x_{i}}%
-\dot{x_{i+1}})+\kappa_{i-1}(x_{i}-x_{i-1})+\kappa_{i}(x_{i}-x_{i+1}%
)=0,\;2\leq i\leq n-1\;,\nonumber\\
M_{n}\ddot{x}_{n}+b_{n-1}(\dot{x}_{n}-\dot{x}_{n-1})+b_{n}(\dot{x_{n}}%
-\dot{x_{nes}})+\kappa_{n-1}(x_{n}-x_{n-1})+\kappa_{nes}(x_{n}-x_{nes}%
)^{3}=0,\\
M_{nes}\ddot{x}_{nes}+b_{n}(\dot{x}_{nes}-\dot{x}_{n})+\kappa_{nes}%
(x_{nes}-x_{n})^{3}=0\;.
\end{gather}\label{eq_ndof01}
Defining a non-dimensional time,%

\end{subequations}
\begin{equation}
\tau=\sqrt{\frac{\kappa_{0}}{M_{1}}}t. \label{eq_ndof02}%
\end{equation}
Equations (\ref{eq_ndof01} (a-d)) can now be written in a non-dimensional
form as
\begin{subequations}
\begin{gather}
x_{1}^{\prime\prime}+2\zeta_{0}x_{1}^{\prime}+2\zeta_{1}(x_{1}^{\prime}%
-x_{2}^{\prime})+x_{1}+k_{1}(x_{1}-x_{2})=0,\\
\mu_{i}x_{i}^{\prime\prime}+2\zeta_{i-1}(x_{i}^{\prime}-x_{i-1}^{\prime
})+2\zeta_{i}(x_{i}^{\prime}-x_{i+1}^{\prime})+k_{i-1}(x_{i}-x_{i-1}%
)+k_{i}(x_{i}-x_{i+1})=0,\;2\leq i\leq n-1,\nonumber\\
\mu_{n}x_{n}^{\prime\prime}+2\zeta_{n-1}(x_{n}^{\prime}-x_{n-1}^{\prime
})+2\zeta_{n}(x_{n}^{\prime}-x_{nes}^{\prime})+k_{n-1}(x_{n}-x_{n-1}%
)+C(x_{n}-x_{nes})^{3}=0,\\
\epsilon x_{nes}^{\prime\prime}+2\zeta_{n}(x_{nes}^{\prime}-x_{n}^{\prime
})+C(x_{nes}-x_{n})^{3}=0,
\end{gather}\label{eq_ndof03}
where, a $^{\prime}$ denotes a derivative with respect to $\tau$ and $\mu
_{i}=\frac{M_{i}}{M_{1}}$ ($1\leq i\leq n$), $\epsilon=\frac{M_{nes}}{M_{1}}$,
$\zeta_{i}=\frac{b_{i}}{2\sqrt{M_{1}\kappa_{0}}}$ ($0\leq i\leq n$),
$k_{i}=\frac{\kappa_{i}}{\kappa_{0}}$ ($1\leq i\leq n-1$) and $C=\frac
{\kappa_{nes}}{\kappa_{0}}$.

\subsection*{Analysis Using Complexification-Averaging}

In a process similar to the one followed in the earlier sections, the
non-dimensional equations of motion given by Eq. (\ref{eq_ndof03}) can be
further analyzed using the technique of complexification-averaging. For this
purpose, new variables are defined as%

\end{subequations}
\begin{subequations}
\begin{gather}
\psi_{i}=x_{i}^{\prime}(\tau)+\omega jx_{i}(\tau),\;1\leq i\leq n,\\
\psi_{nes}=x_{nes}^{\prime}(\tau)+\omega jx_{nes}(\tau),
\end{gather}\label{eq_ndof04}
and further substituting%

\end{subequations}
\begin{subequations}
\begin{gather}
\psi_{i}=\phi_{i}e^{j\omega\tau}\;,\;1\leq i\leq n,\\
\psi_{nes}=\phi_{nes}e^{j\omega\tau}.
\end{gather}\label{eq_ndof05}
Substituting Eqs. (\ref{eq_ndof04}) and (\ref{eq_ndof05}) into Eq.
(\ref{eq_ndof03}) and averaging over the fast time scale $\tau$, the following
equations are obtained:%

\end{subequations}
\begin{subequations}
\begin{gather}
\phi_{1}^{\prime}+\left(  \zeta_{0}+\frac{j\omega}{2}-\frac{j}{2\omega
}\right)  \phi_{1}+\zeta_{1}(\phi_{1}-\phi_{2})-\frac{jk_{1}}{2\omega}%
(\phi_{1}-\phi_{2})=0,\\
\mu\phi_{i}^{\prime}+\zeta_{i-1}(\phi_{i}-\phi_{i-1})+\zeta_{i}(\phi_{i}%
-\phi_{i+1})+\frac{j\omega}{2}\mu_{i}\phi_{i}-\frac{jk_{i-1}}{2\omega}%
(\phi_{i}-\phi_{i-1})-\frac{jk_{i}}{2\omega}(\phi_{i}-\phi_{i+1})=0,\;2\leq
i\leq n-1,\nonumber\\
\mu_{n}\phi_{n}^{\prime}+\zeta_{n-1}(\phi_{n}-\phi_{n-1})+\zeta_{n}(\phi
_{n}-\phi_{nes})+\frac{j\omega}{2}\mu_{n}\phi_{n}-\frac{jk_{n-1}}{2\omega
}(\phi_{n}-\phi_{n-1})-\frac{3jC}{8\omega^{3}}|\phi_{n}-\phi_{nes}|^{2}%
(\phi_{n}-\phi_{nes})=0,\\
\epsilon\phi_{nes}^{\prime}+\zeta_{n}(\phi_{nes}-\phi_{n})+\frac{j\omega}%
{2}\epsilon\phi_{nes}-\frac{3jC}{8\omega^{3}}|\phi_{nes}-\phi_{n}|^{2}%
(\phi_{nes}-\phi_{n})=0.
\end{gather}\label{eq_ndof06}
Based on Eqs. (\ref{eq_ndof06}), some new variables can now be defined as%

\end{subequations}
\begin{subequations}
\begin{gather}
u_{i}=\phi_{i}-\phi_{i+1},\;1\leq i\leq n-1,\\
u_{n}=\phi_{n}-\phi_{nes},\\
u_{n+1}=\phi_{1}+\sum_{i=2}^{n}\mu_{i}\phi_{i}+\epsilon\phi_{nes}.
\end{gather}\label{eq_ndof07}

There are now $n+1$ variables, of which $u_{1}$ to $u_{n}$ represent the
relative displacements between the masses and $u_{n+1}$ is the motion of the
center of mass the system. The main variable of interest here is $u_{n}$ as it
represents the relative displacement of the primary system and the NES. It is
an approximate measure of the energy being dissipated or removed from the
primary system. Using the definitions given in Eqs. (\ref{eq_ndof07}) and
using Eqs. (\ref{eq_ndof06}), the equations of motion for the variables
$u_{i}$ ($i = 1 \;$ to $\; n+1$) can be written as,%

\end{subequations}
\begin{subequations}
\begin{gather}
u_{i}^{\prime}+\sum_{j=1}^{n+1}c_{ij}u_{j}=0,\;1\leq i\leq n-2,\\
u_{n-1}^{\prime}+\sum_{j=1}^{n+1}c_{n-1j}u_{j}+\frac{3jC}{8\mu_{n}\omega^{3}%
}|u_{n}|^{2}u_{n}=0,\\
u_{n}^{\prime}+\sum_{j=1}^{n+1}c_{nj}u_{j}-\frac{3jC(\mu_{n}+\epsilon)}%
{8\mu_{n}\epsilon\omega^{3}}|u_{n}|^{2}u_{n}=0,\\
u_{n+1}^{\prime}+\sum_{j=1}^{n+1}c_{n+1j}u_{j}=0,
\end{gather}\label{eq_ndof08}
where $c_{ij}$, ($1\leq i,j\leq n+1$), are constant coefficients. Assuming
that a starting impulse velocity of $v_{0}$ is given to the first mass, the
only two variables having non-zero initial conditions are,%

\end{subequations}
\begin{subequations}
\begin{gather}
u_{1}(0)=\phi_{1}(0)-\phi_{2}(0)=v_{0},\\
u_{n+1}(0)=\phi_{1}(0)+\sum_{i=2}^{n}\mu_{i}\phi_{i}(0)+\epsilon\phi
_{nes}(0)=v_{0}.
\end{gather}\label{eq_ndof09}

It can be observed from Eq. (\ref{eq_ndof08}) that there are $n-2$
linear differential equations and 2 nonlinear differential equations (for
$u_{n-1}$ and $u_{n}$). Out of the nonlinear equations, the variable $u_{n}$
is of the most interest since it is representative of effectiveness of the NES. The system of equations given in Eq.
(\ref{eq_ndof08}) are also referred to as the slow-flow equations of motion as
they are derived for a slow-time scale.

\subsection*{Solution of Slow-Flow Equations of Motion}

Considering the coefficients $c_{ij}$ of Eq. (\ref{eq_ndof08}), it is observed
that $c_{in}$, $(1\leq i\leq n-2\;,$ and $,\;i=n+1\;)$ are of $O(\epsilon)$ .
Therefore, the the slow flow equations given in Eq. (\ref{eq_ndof08}) can now
be approximated as%

\end{subequations}
\begin{subequations}
\begin{gather}
u_{n-1}^{\prime}+c_{n-1n-1}u_{n-1}+c_{n-1n}u_{n}+\bar{d_{2}}\bar{u}_{p}%
+\frac{3jC}{8\mu_{n}\omega^{3}}|u_{n}|^{2}u_{n}=0,\\
u_{n}^{\prime}+c_{nn-1}u_{n-1}+c_{nn}u_{n}+\bar{d_{3}}\bar{u}_{p}%
-\frac{3jC(\mu_{n}+\epsilon)}{8\mu_{n}\epsilon\omega^{3}}|u_{n}|^{2}u_{n}=0,\\
\bar{u}_{p}^{\prime}+\bar{A}\bar{u}_{p}+\bar{d_{1}}u_{n-1}+O(\epsilon)=0,
\end{gather}\label{eq_ndof10}
where, $\bar{u}_{p}$, is a column vector of $u_{i}$, $(1\leq i\leq n-2)$ and
$i=n+1$. Thus, the linear equations have been collected under a single vector
variable, $\bar{u}_{p}$. The variable $\bar{A}$ is a $(n-1\times n-1)$ matrix
containing the relevant coefficients and $\bar{d_{3}}$, $\bar{d_{2}}$ and
$\bar{d_{3}}$ are $(1\times n-1)$ row vectors. The Eqs. (\ref{eq_ndof10}) are
quite similar to Eqs. (\ref{eq_2dof15}) which have been extensively studied in
this work (the only difference is that Eq. (\ref{eq_ndof10})(c) is a linear
vector equation rather than being a linear scalar equation). Thus, the
analysis done for the two-degree-of-freedom system can be used for the N-degree-of-freedom system without
much additional effort.
\end{subequations}
\bibliographystyle{elsarticle-num-names}

\end{document}